\documentclass{article}

\usepackage{arxiv}

\usepackage[utf8]{inputenc} % allow utf-8 input
\usepackage[T1]{fontenc}    % use 8-bit T1 fonts
\usepackage{hyperref}       % hyperlinks
\usepackage{url}            % simple URL typesetting
\usepackage{booktabs}       % professional-quality tables
\usepackage{amsfonts}       % blackboard math symbols
\usepackage{nicefrac}       % compact symbols for 1/2, etc.
\usepackage{microtype}      % microtypography
\usepackage{lipsum}
\usepackage{graphicx}

\usepackage{amsmath}
\usepackage{amssymb}
\usepackage{amsthm}

\usepackage{caption}

\usepackage{graphicx}
\usepackage{wrapfig}
\usepackage{subfigure}

\usepackage{substr}
\usepackage[tight]{minitoc}

\usepackage[all]{xy}

\usepackage{algorithm}
\usepackage{algpseudocode}
\usepackage{xcolor}
\usepackage{enumitem}

\def\R{\mathbb R}
\def \E {\mathbb E}
\newcommand{\circlednum}[1]{\textcircled{#1}}

\newtheorem{teo}{Theorem}

\theoremstyle{remark}
\newtheorem{remark}{Remark}

\newtheorem{assum}{Assumption}

\algnewcommand{\Inputs}[1]{%
  \State \textbf{Inputs:}
  \Statex \hspace*{\algorithmicindent}\parbox[t]{.8\linewidth}{\raggedright #1}
}
\algnewcommand{\Initialize}[1]{%
  \State \textbf{Initialize:}
  \Statex \hspace*{\algorithmicindent}\parbox[t]{.8\linewidth}{\raggedright #1}
}
\algdef{SE}[DOWHILE]{Do}{doWhile}{\algorithmicdo}[1]{\algorithmicwhile\ #1}

\graphicspath{ {./images/} }

\title{Accelerated zero-order SGD under high-order
smoothness and overparameterized regime}

\author{
 Georgii Bychkov \\
  Lomonosov MSU, ISP RAS\\
  Moscow, Russia \\
  \texttt{georgy.bychkov@graphics.cs.msu.ru} \\
   \And
 Darina Dvinskikh \\
  HSE University, MIPT, ISP RAS\\
  Moscow, Russia \\
  \texttt{dmdvinskikh@hse.ru} \\
  \And
 Anastasia Antsiferova \\
  MSU Institute for Artificial Intelligence, ISP RAS\\
  Moscow, Russia \\
  \texttt{aantsiferova@graphics.cs.msu.ru} \\
  \And
 Alexander Gasnikov \\
  Innopolis University, MIPT, IITP RAS\\
  Innopolis, Republic of Tatarstan, Russia \\
  \texttt{gasnikov@yandex.ru} \\
  \And
 Aleksandr Lobanov \\
  MIPT, Skoltech, ISP RAS\\
  Dolgoprudny, Moscow Oblast,  Russia \\
  \texttt{lobbsasha@mail.ru} \\
}

\begin{document}
\maketitle
\begin{abstract}
We present a novel gradient-free algorithm to solve a convex stochastic optimization problem, such as those encountered in medicine, physics, and machine learning (e.g., adversarial multi-armed bandit problem), where the objective function can only be computed through numerical simulation, either as the result of a real experiment or as feedback given by the function evaluations from an adversary. Thus we suppose that only a black-box access to the function values of the objective is available, possibly corrupted by adversarial noise: deterministic or stochastic. The noisy setup can arise naturally from modeling randomness within a simulation or by computer discretization, or when exact values of function are forbidden due to privacy issues, or when solving non-convex problems as convex ones with an inexact function oracle. By exploiting higher-order smoothness, fulfilled, e.g., in logistic regression, we improve the performance of zero-order methods developed under the assumption of classical smoothness (or having a Lipschitz gradient). The proposed algorithm enjoys optimal oracle complexity and is designed under an overparameterization setup, i.e., when the number of model parameters is much larger than the size of the training dataset. Overparametrized models fit to the training data perfectly while also having good generalization and outperforming underparameterized models on unseen data. We provide convergence guarantees for the proposed algorithm under both types of noise. Moreover, we estimate the maximum permissible adversarial noise level that maintains the desired accuracy in the Euclidean setup, and then we extend our results to a non-Euclidean setup. Our theoretical results are verified on the logistic regression problem.
\end{abstract}

% keywords can be removed

\section{Introduction}
\label{sec:introduction}
\setcounter{equation}{0}
In this paper we focus on the following stochastic convex optimization problem
\begin{equation}\label{eq:main_prob}
\min\limits_{x\in \R^d}\{f(x) := \E_\xi f(x, \xi) \}, 
\end{equation}
where $f$ is convex in $x$, and $\E_\xi$ is the expectation,  with respect to $\xi$ with unknown distribution. 

Motivated by   different applications where gradient information  is unavailable or calculating gradients is prohibitively expensive, we seek to solve \eqref{eq:main_prob}
 only by accessing noisy evaluations of $f(x, \xi)$  \cite{polyak1987introduction,granichin2003randomized,risteski2016algorithms,vasin2021stopping}, possibly adversarial noisy.  
Methods that query only function values for given inputs, rather than gradients, are known as zero-order, or gradient-free methods, \cite{polyak1987introduction,spall2003introduction,conn2009introduction,duchi2015optimal,shamir2017optimal,nesterov2017random,gasnikov2017stochastic,beznosikov2020gradient,gasnikov2022power}  and are actively studied due to the vast number of applications.
Indeed,  in medicine, biology, and physics, the objective function can only be computed through numerical simulation or as the result of a real experiment, i.e., automatic differentiation cannot  be used to compute function derivatives. 
 We also mention classical problem of adversarial multi-armed bandit \cite{flaxman2004online,bartlett2008high,bubeck2012regret}, where a learner receives a feedback  given by the function evaluations from an adversary.    
 Usually, a black-box function we are optimizing, is affected by stochastic or computational noise. This noise can arise naturally from modeling randomness within a simulation or by computer discretization. The Noisy setup is also useful in cases where exact values of function are forbidden due to privacy issues. Nowadays, there are a lot of adversarial attacks that influence oracle output. An attacker can implement a target function with backdoor attacks~\cite{ge2021anti}. 
 Thus, we address two concepts of adversarial noises: \emph{adversarial deterministic noise} \cite{dvinskikh2022gradient} and \emph{adversarial stochastic noise}  \cite{lobanov2023stochastic}. Adversarial stochastic noise makes the gradient approximation unbiased, so it will only accumulate in the variance. Adversarial deterministic noise, on the other hand, accumulates in both bias and variance. Adversarial deterministic noise has a stronger effect on the optimization process than adversarial stochastic noise.

%One of the challenges in zero-order optimization is to find the optimal trade-off between the number of oracle calls and the error of the solution.
%This paper proposes a method for solving a zero-order stochastic optimization task by considering overparametrized models and convex high-smooth problems. 

There is a large amount  of literature on zero-order optimization, to name a few we mention \cite{duchi2015optimal} for the smooth case (differential and with
Lipschitz-continuous gradients), for the non-smooth case we refer to \cite{shamir2017optimal}. These papers  uses the exact function evaluations that can be infeasible in some applications. Moreover,  non-convex problems are often considered as convex ones   with an inexact function  oracle. Zero-order methods with noisy function evaluations are studied in \cite{polyak1987introduction,granichin2003randomized,bayandina2018gradient,beznosikov2020gradient,risteski2016algorithms,vasin2021stopping}.
%focused on more common settings: underparametrized models [], non-smooth [] and smooth~\cite{lobanov2023accelerated} functions, oracle without noise [], etc. 
In contrast to these works we consider 
overparametrized regime  \cite{jacot2018neural,  allen2019learning, belkin2019reconciling}, i.e., when the number of model parameters is much larger than the size of the training dataset. 
Overparametrized models  fit to the training data perfectly while also having good generalization and outperforming underparameterized models on unseen data.
%This is known in machine learning community as double descent phenomenon. 
For instance, as stated in State of AI report, DeepMind researchers found that language models are significantly undertrained and there is a lack of labeled training data~\cite{hoffmann2022training, eldan2023tinystories}. 
Overparametrized models were studied in \cite{lobanov2023accelerated} for smooth and convex functions. 
 In contrast to \cite{lobanov2023accelerated}, we  consider  higher-order smoothness  \cite{polyak1990optimal, bach2016highly, akhavan2023gradient},  fulfilled e.g.,  in logistic regression, to improve the performance of zero-order  methods. To tackle higher-order smoothness, we use kernel-based gradient approximation \cite{bach2016highly}.
%We consider the more narrow setting for several reasons. %Firstly, overparametrized models frequently arise in machine learning due to the limited amount of labeled data. 
%Secondly, non-convex tasks can be reformulated to convex, which makes the solution more efficient []. 
We summarize our contributions as follows:
\begin{itemize}\itemsep=-2pt
\item We propose a novel gradient-free algorithm enjoying optimal oracle complexity under certain batch sizes.
\item  We  estimate the maximum permissible adversarial noise level that maintains the desired accuracy in the Euclidean setup, and then  we extend our results to a non-Euclidean setup.
\item We confirm our theoretical results on machine learning problem. The numerical experiments demonstrate the superiority of the proposed method over AZO-SGD \cite{lobanov2023accelerated}, designed for  smooth objective functions.
\end{itemize}

%  Imprecise function values can be returned if the exact calculation is costly, and we want to reduce the computational complexity of optimization. Also, the oracle cannot return exact values due to privacy issues. Another big area covered by stochastic optimization is the optimization in case of adversarial attacks that influence oracle output. An attacker can implement a target function with backdoor attacks~\cite{bib:backdoor}.
%Finally, we consider stochastic optimization with an oracle that returns a response containing some noise bounded by $\Delta$. Imprecise function values can be returned if the exact calculation is costly, and we want to reduce the computational complexity of optimization. Also, the oracle cannot return exact values due to privacy issues. Another big area covered by stochastic optimization is the optimization in case of adversarial attacks that influence oracle output. An attacker can implement a target function with backdoor attacks~\cite{bib:backdoor}.
%This paper also focuses on finding the maximum admissible noise level at which the convergence to the desired accuracy can still be guaranteed. 
%\subsection{Contributions}

\paragraph{Paper organization.}

The structure of the paper is the following. In Section \ref{sec:definition}, we provide background material, notation, and assumptions. In Section \ref{sec:main}, we  present the main algorithm to solve  \eqref{eq:main_prob} along with the convergence guarantees. Section \ref{sec:experiments} provides the numerical experiments demonstrating the superiority of our novel algorithm
on the logistic regression. In Section \ref{sec:conclusion}, we summarize the paper and formulate conclusions.

\section{Preliminaries}\label{sec:definition}
In this section, we provide all the notation and assumptions necessary to prove the convergence of our novel algorithm.
\paragraph{Notation.}

We use $\|x\|_p$ to denote the $\ell_p$-norm, for the Euclidean norm we use the shorthand notation $\|x\|_2 = \|x\|$. We use $\langle x, y \rangle := \sum_{k=1}^d x_iy_i$ as the standard inner product of $x, y \in \mathbb{R}^d$, where $x_i$ and $y_i$ are the $i$-th component of $x$ and $y$ respectively. We denote by $S_p^d(r) := \{x\in \mathbb{R}^d :\|x\|_p = r\}$ and $B_p^d(r) := \{x\in \mathbb{R}^d :\|x\|_p \le r\}$ the sphere and the ball in the $\ell_p$-norm respectively. Symbol $\lesssim$ stands for the  asymptotic inequality.

%\paragraph{Assumptions.}
%Our paper focuses on the case when we obtain an inexact gradient value when calling the oracle, which is referred to as the biased oracle~\cite{lobanov2023accelerated}. It is encountered in adversarial attacks or as a gradient approximation when the exact gradient is difficult.

%, where $K(u)$ is defined in Section \ref{sec:algorithm}.

%\subsection{Assumptions}
% Same as paper Accelerated Zero-Order SGD Method for Solving the Black Box Optimization Problem under “Overparametrization” Condition

%\begin{assum}[Convexity]
%\label{assum:conv}
%Function $f$ is convex if
%\begin{equation}
%\forall x, y, \xi \quad f(y, \xi) \ge f(x, \xi) + \langle \nabla f(x, \xi), y - x\rangle.
%\end{equation}
%\end{assum}

\begin{assum}[$L$-smoothness]\label{ass:L-smoothness}
 $f(x,\xi)$ is $L$-Lipschitz smooth, or has  $L$-Lipschitz continuous gradient, if $f(x,\xi)$ is continuously differentiable, with respect to $x$, and its gradient satisfies Lipschitz condition for any $\xi$ and $x,y \in \R^d$
\begin{equation*}
\|\nabla f(x,\xi) - \nabla f(y,\xi) \|\le L \|x-y\|.
\end{equation*}
\end{assum}
%\begin{assum}[$L$-smoothness]
%\label{assum:smooth}
%Function $f$ is $L$-smooth if
%\begin{equation}
%\forall x, y, \xi \quad \exists L: \|f(x, \xi) - f(y, \xi)\| \le L\|x - y\|.
%\end{equation}
%\end{assum}
%\begin{assum}
%\label{assum:min}
%The function $f(x)$ is a convex and has the minimum value $f^* = \min_x f(x)$, which is attained at a point $x^*$ with $\|x^* - x_0\| \le R$, where $x_0$ is a starting point.
%\end{assum}
%These assumptions are commonly used in various optimization papers (like Lobanov et al.\cite{bib:LobanovAccelerated}, Woodworth et al.\cite{bib:WoodworthEven}, Ilandarideva et al.\cite{bib:IlandaridevaAccelerated}).
The next assumption is an  extension of Assumption \ref{ass:L-smoothness}.
\begin{assum}[Higher-order smoothness]\label{ass:High-smoothness}
Let $l$ denote the maximal integer number strictly less than $\beta$. Let $\mathcal{F}_\beta(L)$ denote the set of all functions $f : R^d \to R$ which are differentiable $l$ times and for any $\xi$ and $ x,y \in \R^d$ the Hölder-type condition: 
\begin{equation*}
\left|f(x,\xi) - \sum_{0 \le |n| \le l}\frac{1}{n!}D^nf(y,\xi)(x - y)^n\right| \le L_\beta \|x - y\|,
\end{equation*}
where $L_\beta > 0$, the sum is over multi-index $n = (n_1,...,n_d) \in \mathbb{N}^d$, the notation $n! = n_1!...n_d!$, and 
\begin{equation*}
D^nf(y,\xi) v^n = \frac{\partial^{|n|}f(y,\xi)}{\partial^{n_1}x_1...\partial^{n_d}x_d}v_1^{n_1}...v_d^{n_d},
\end{equation*}
where $|n| = n_1+...+n_d$, for all $ v =(v_1,...,v_d) \in \mathbb{R}^d$.
\end{assum}

When $\beta = 2$, Assumption \ref{ass:High-smoothness} turns into Assumption \ref{ass:L-smoothness}.

%\subsection{Assumptions on the gradient oracle}
% Same as paper Accelerated Zero-Order SGD Method for Solving the Black Box Optimization Problem under “Overparametrization” Condition

%\begin{fed}[Gradient oracle] 
%\label{fed:biased}
%A map $\mathbf{g}:\mathbb{R}^d \times \mathcal{D} \to \mathbb{R}^d$ s.t.
%\begin{equation*}
%\mathbf{g}(x, \xi) = \nabla f(x, \xi) + \mathbf{b}(x),
%\end{equation*}
%where $\mathbf{b}: \mathbb{R}^d \to \mathbb{R}^d$ such that for all $ x \in \mathbb{R}^d: \|\mathbf{b}(x)\|^2 \le \zeta^2.$
%\end{fed}

The next assumption means that in overparametrized regime  the variance of the stochastic gradient estimates at optimal  point $x^* = \arg\min\limits_{x\in \R^d}f(x)$ can be upper bounded by using the minimum value $ f(x^*)$ . 
\begin{assum}[\cite{woodworth2021even}, Lemma 3]\label{assum:overparam}
There exists $\sigma_*^2 \ge 0$ such  that
\begin{equation*}
\mathbb{E}[\|\nabla f(x^*, \xi) - \nabla f(x^*)\|^2] \le \sigma_*^2.
\end{equation*}
\end{assum}

We will use this assumption for the gradient estimates.

%\subsection{Adversarial noise concept}
% Write about deterministic and stochastic adversarial noise and their differences

%We estimate $\Delta$ for different noise variants.

%\section{Problem Statement}  \label{sec:definition}

%\begin{teo}[Convergence of Biased AC-SA \cite{bib:LobanovAccelerated}]
%Let $f$ satisfy Assumptions
%\ref{assum:conv},\ref{assum:smooth},\ref{assum:min} and gradient oracle from Definition \ref{fed:biased} satisfy Assumption \ref{assum:overparam}, then Biased
%AC-SA algorithm guarantees convergence with a universal constant $c$
%\begin{equation}
%\mathbb{E}[f(x_N^{ag}) - f^*] \le c\left( \frac{LR^2}{N^2} + \frac{LR^2}{BN} + \frac{\zeta R}{\sqrt{BN}} + b R + %\frac{b^2}{2L}N\right)
%\end{equation}
%\end{teo}

\section{Main results: the AZO-SGD-HS algorithm}
\label{sec:main}
%\subsection{}
%\label{sec:algorithm}
% Write modification on gradient estimation on AZO-SGD

In this section, we present a novel algorithm, named AZO-SGD-HS (see Algorithm \ref{alg:azo-hs}) and provide its theoretical guarantees of convergence when the adversarial noise is deterministic and when it is stochastic. We also report on the maximum permissible adversarial noise level  $\Delta$. 
%The AZO-SGD-HS algorithm generalizes the AZO-SGD \cite{lobanov2023accelerated}, designed for smooth  gradient-free overparametrized setup.  
%Our algorithm, as well the AZO-SGD, is based on the first-order algorithm \cite{woodworth2021even} and uses  the two-point gradient approximation. The main difference of our algorithm from the AZO-SGD is that   it uses kernel-based gradient approximation \cite{bach2016highly}, which is more appropriate  for higher-order smoothness problems.

\subsection{Deterministic adversarial noise}
% Write theorem with proof in the case of deterministic adversarial noise

Let us suppose that   a zero-order oracle is corrupted by a deterministic  noise $\delta(x)$, bounded by a some constant  $\Delta: |\delta(x)| \le \Delta$  for all $x\in \R^n$ 
\begin{equation}\label{eq:zero_order_deter}
f_\delta(x, \xi) = f(x, \xi) + \delta(x).
\end{equation}
Then by using two-point noisy zero-order oracle \eqref{eq:zero_order_deter} the gradient can be estimated as follows 
\begin{equation}\label{eq:grad_approx}
\mathbf{g}(x, \mathbf{e}, r, \xi) = \frac{d}{2h}(f_\delta(x + hr\mathbf{e}, \xi) - f_\delta(x - hr\mathbf{e}, \xi))K(r)\mathbf{e},
\end{equation}
where $h > 0$ is a smoothing parameter, $\mathbf{e} \in S_2^d(1)$ is a vector uniformly distributed on the Euclidean unit sphere,  $K:[-1, 1] \to \mathbb{R}$ is a kernel function that satisfies
\[
E[K(u)] = 0, \quad E[uK(u)] = 1,\quad  E[|u|^\beta|K(u)|] < \infty \text{ and for all } j = 2, ..., l \quad  E[u_
jK(u)] = 0.
\]

\begin{algorithm}{}
	\caption{Accelerated zero-order SGD for high-order smoothness (AZO-SGD-HS)}
	\label{alg:azo-hs}
	
        \begin{algorithmic}[1]
        \Inputs{Starting points $x_0^{ag} = x_0 = \vec{0} \in \mathbb{R}^d$,  number of iterations $N$, batch size $B$, smoothing parameter $h > 0$, kernel  $K(\cdot)$.}
        \For{$k$ = $0$ to $N-1$}
        \State $\beta_k = 1 + \frac{k}{6}$; $\gamma_k = \gamma(k+1)$ where $\gamma=\min\left\{\frac{1}{12L}, ~ \frac{B}{24L(N+1)},~ \sqrt{\frac{BR^2}{Lf^*N^3}}\right\}$.
        \State $x_k^{md} = \beta_k^{-1}x_k + (1 - \beta_k^{-1})x_k^{ag}$.
        \State Sample $\mathbf{e_i} \in S_2^d(1)$, $r_i = {\rm Uniform}[-1, 1]$ and $\xi_i$ for $i=1, ..., B$. 
        \State Compute batched gradient approximation  
        \begin{equation}\label{eq_batched_grad}
            g_k =\frac{1}{B}\sum_{i=1}^B 
        \mathbf{g}(x_k^{md},  \mathbf{e_i}, r_i, \xi),
        \end{equation}
            \hspace{1.3em} where $\mathbf{g}$ is defined in 
        \begin{itemize}[left=1.3em]
            \item deterministic noise:  \eqref{eq:grad_approx}
            \item stochastic noise: \eqref{eq:grad_stoch_noise}
        \end{itemize}
 
        \State $\overline{x}_{k+1} = x_k - \gamma g_k.$
        \State $x_{k+1} = \min\left\{1, \frac{R}{\|\overline{x}_{k+1}\|}\right\}\overline{x}_{k+1}.$
        \State $x_{k+1}^{ag} = \beta_k^{-1}x_{k+1} + (1 - \beta_k^{-1})x_k^{ag}.$
        \EndFor
        
    \State \Return $x_N^{ag}.$
      \end{algorithmic}
\end{algorithm}

The next theorem states the convergence guarantee for AZO-SGD-HS with  two-point zero-order oracle  \eqref{eq:zero_order_deter}. We
 will use the shorthand notations $\kappa_\beta = \int |u|^\beta |K(u)|du$ and $\kappa = \int |K(u)^2|du$.
\begin{teo}
\label{teo:convergence_deterministic}
Let \(f(\cdot,\cdot)\) satisfy Assumption \ref{ass:High-smoothness} with parameter $\beta$. Let smoothing parameter be $h \le \left({\varepsilon}/{(\kappa_{\beta}LR})\right)^{{1}/{(\beta - 1)}}$ then    Assumption \ref{assum:overparam} is satisfied with gradient approximation  \eqref{eq_batched_grad}. Let the maximum admissible level of adversarial deterministic noise be 
\begin{equation*}
\Delta \leq \frac{\varepsilon^{1 + {1}/{(\beta - 1)}}}{d\kappa_{\beta}^{{1}/{(\beta - 1)}}R^{1 + {1}/{(\beta - 1)}} L^{{1}/{(\beta - 1)}}}.
\end{equation*}
Let $x_N^{ag}$ be the output of the AZO-SGD-HS then 
\[\mathbb{E}[f(x_N^{ag}) - f(x^*)] \le \varepsilon\] 
in at most the following number of iterations $N$ and  oracle calls $T$, where $R$ will denote a radius of a solution $x^*$: 
$\|x^* - x_0\| \le R$, and $x_0$ is a starting point:

 \begin{itemize}
     \item $B = 1$
\begin{equation*}
N = \mathcal{O}\left(\max\left(\frac{LR^2}{\varepsilon}, \frac{d\kappa\sigma_*^2R^2}{\varepsilon^2}\right)\right),\ T = \mathcal{O}\left(\max\left(\frac{LR^2}{\varepsilon}, \frac{d\kappa\sigma_*^2R^2}{\varepsilon^2}\right)\right)
\end{equation*}

     \item 
     $N > B > 1$
\begin{equation*}
N = \mathcal{O}\left(\max\left(\frac{LR^2}{B\varepsilon}, \sqrt{\frac{LR^2}{\varepsilon}}, \frac{d\kappa\sigma_*^2R^2}{B\varepsilon^2}\right)\right)
\end{equation*}
\begin{equation*}
T = \mathcal{O}\left(\max\left(\frac{LR^2}{\varepsilon}, B\sqrt{\frac{LR^2}{\varepsilon}}, \frac{d\kappa\sigma_*^2R^2}{\varepsilon^2}\right)\right) 
\end{equation*}

     \item 

$B = N$
\begin{equation*}
N = \mathcal{O}\left(\max \left(\sqrt{\frac{LR^2}{\varepsilon}}, \frac{\sqrt{d\kappa}\sigma_*R}{\varepsilon}\right)\right),\ T = \mathcal{O}\left(\max \left( \frac{LR^2}{\varepsilon}, \frac{d\kappa\sigma_*^2R^2}{\varepsilon^2}\right)\right)  
\end{equation*}

     \item 
$B > N$
\begin{equation*}
N = \mathcal{O}\left( \sqrt{\frac{LR^2}{\varepsilon}}\right) ,\ T = \mathcal{O}\left(\max \left(\frac{LR^2}{\varepsilon} , \frac{d\kappa\sigma_*^2R^2}{\varepsilon^{2}} , \frac{\kappa_{\beta}^{2/(\beta - 1)}\kappa d^2L^{2/(\beta - 1)}\Delta^2R^{2 + 2/(\beta - 1)}}{\varepsilon^{2 + 2/(\beta - 1)}}\right) \right).
\end{equation*}
     
 \end{itemize}
 
\end{teo}

\textit{Sketch of the proof. }
To find the upper bounds for the maximum admissible level of adversarial noise $\Delta$, smoothing parameter $h$, the asymptotic for the number of iteration $N$, and the total number of oracle calls $T$.
The proof consists in  upper bounding  the bias and second moment of the gradient approximation. 

\textbf{Step 1.} (bias of gradient approximation). 

We calculate the bias of the gradient approximation $\|\mathbb{E}[\mathbf{g}(x_k, \mathbf{e}, r, \xi)] - \nabla f(x_k)\|$.
\begin{align*}
 b &= \|\mathbb{E}[\mathbf{g}(x_k, \mathbf{e}, r, \xi)] - \nabla f(x_k)\| = \left\|\mathbb{E}\left[\frac{d}{2h} (f_\delta(x_k + hr\mathbf{e}, \xi) - f_\delta(x_k - hr\mathbf{e}, \xi)) K(r) \mathbf{e}\right] - \nabla f(x_k)\right\|  \\
& \stackrel{\text{\circlednum{1}}}{=}  \left\|\mathbb{E}\left[ \frac{d}{2h}(f(x_k + hr\mathbf{e}, \xi))K(r) \mathbf{e} + \delta(x_k + h \mathbf{e})\right] - \nabla f(x_k)\right\| \\
& \stackrel{\text{\circlednum{2}}}{\le} \left\|\mathbb{E}\left[\frac{d}{2h}(f(x_k + hr\mathbf{e}, \xi))K(r) \mathbf{e}\right] - \nabla f(x_k)\right\| + \frac{d\Delta}{h}\\
& \stackrel{\text{\circlednum{3}}}{=}\|\mathbb{E}([\nabla f(x_k + h r \mathbf{u}, \xi)r K(r)] - \nabla f(x_k))\| + \frac{d\Delta}{h}\\
& \stackrel{\text{\circlednum{4}}}{=}\sup_{z \in S_2^d(1)}\mathbb{E}[(\nabla_z f(x_k + h r \mathbf{u}, \xi) - \nabla_z f(x_k))r K(r)] + \frac{d\Delta}{h} \\
& \stackrel{\text{\circlednum{5}}}{\le} \kappa_{\beta}h^{\beta - 1} \frac{L}{(l-1)!}\mathbb{E}[\|\mathbf{u}\|^{\beta - 1}] + \frac{d\Delta}{h} \le \kappa_{\beta}h^{\beta - 1} \frac{L}{(l-1)!}\frac{d}{d + \beta - 1} + \frac{d\Delta}{h}\lesssim \kappa_{\beta}L h^{\beta - 1} + \frac{d\Delta}{h},
\end{align*}
where $\mathbf{u} \in B_2^d(1)$, \circlednum{1}: distribution of $\mathbf{e}$ is symmetric; \circlednum{2}: triangle inequality and bounded noise $|\delta(x)| \le \Delta$; \circlednum{3}: a version of the Stokes’
theorem \cite{zorich2016mathematical} (see Section 13.3.5, Exercise 14a); \circlednum{4}: norm of the gradient is the supremum of directional derivatives $\nabla_z f(x) = \lim_{\varepsilon \to 0}\frac{f(x + \varepsilon z) - f(x)}{\varepsilon}$; \circlednum{5}: Taylor expansion.

\textbf{Step 2.}  (bounding the second moment of gradient approximation).

In this step, we calculate the second moment of the gradient approximation $\mathbb{E}\|\mathbf{g}(x^*, \mathbf{e}, r, \xi)\|^2$.
\begin{align*}
\zeta^2 &= \mathbb{E} [\|\mathbf{g}(x^*, \mathbf{e}, r, \xi)\|^2 ] = \mathbb{E} \left [\left\| \frac{d}{2h}(f_\delta(x^* + hr\mathbf{e}, \xi) - f_\delta(x^* - hr\mathbf{e}, \xi)) K(r) \mathbf{e}\right\|^2\right] \\
& =\frac{d^2}{4h^2}\mathbb{E} [((f_\delta(x^* + hr\mathbf{e}, \xi) - f_\delta(x^* - hr\mathbf{e}, \xi)) K(r))^2]\\
& =\frac{d^2}{4h^2}\mathbb{E} [(f(x^* + hr\mathbf{e}, \xi)- f(x^* - hr\mathbf{e}, \xi) + \delta(x^* + hr\mathbf{e}) - \delta(x^* - hr\mathbf{e}))^2 (K(r))^2] \\
& \stackrel{\text{\circlednum{1}}}{\le} \frac{d^2 \kappa}{2h^2} (\mathbb{E} (f(x^* + hr\mathbf{e}, \xi)- f(x^* - hr\mathbf{e}, \xi))^2 + 2\Delta^2 ) \\
&\stackrel{\text{\circlednum{2}}}{\le} \frac{d^2 \kappa}{2h^2} (\frac{h^2}{d}\mathbb{E} (\|\nabla f(x^* + hr\mathbf{e}, \xi)+\nabla f(x^* - hr\mathbf{e}, \xi)\|^2) + 2\Delta^2)\\
& = \frac{d^2\kappa}{2h^2} (\frac{h^2}{d}\mathbb{E} (\|\nabla f(x^* + hr\mathbf{e}, \xi)+\nabla f(x^* - hr\mathbf{e}, \xi) \pm 2\nabla f(x^*, \xi)\|^2) + 2\Delta^2) \\
& \stackrel{\text{\circlednum{3}}}{\le} 4d\kappa\|\nabla f(x^*, \xi)\|^2 + 4d\kappa L^2h^2\mathbb{E}[\|\mathbf{e}\|^2] + \frac{d^2 \kappa\Delta^2}{h^2}  \\
& \stackrel{\text{\circlednum{4}}}{\le} 4d\kappa\sigma_*^2 + 4d\kappa L^2h^2\mathbb{E}[\|\mathbf{e}\|^2] + \frac{d^2 \kappa\Delta^2}{h^2} = 4d\kappa\sigma_*^2 + 4d\kappa L^2h^2 + \frac{d^2\kappa\Delta^2}{h^2},
\end{align*}
\circlednum{1}: inequality of squared norm of the sum, inequality between positive random variables and independence of the noise; \circlednum{2}: Wirtinger-Poincare inequality; \circlednum{3}: L-smoothness function; \circlednum{4}: overparametrization assumption.

\qed

For the complete proof of this theorem, we refer to Appendix \ref{app:teo_convergence_deterministic}. \\

The Theorem \ref{teo:convergence_deterministic} shows that the iteration complexity is seriously affected by chosen $B$, and a large enough value can make the algorithm have an effective iterative complexity of $\mathcal{O}\left(\sqrt{{LR^2}/{\varepsilon}}\right)$. Total oracle calls are generally worse than AZO-SGD due to the large $B$ needed for optimal iteration complexity. Smoothing parameter value increases with the growth of $\beta$ that is caused by use of kernel approximation that efficiently solves the task. The theorem gives an estimation of the maximum allowable level of adversarial noise. The bound matches the corresponding value from AZO-SGD in the case of $\beta = 2$. Moreover, with the growth of smoothness, the restriction on the noise level becomes weaker, allowing for better-quality solutions.

\begin{remark}
\label{rem:convergence_deterministic}
Theorem \ref{teo:convergence_deterministic} can be generalized to a more general class of problems by using the  $\ell_p$-norm. The number of iterations $N$, the maximum admissible level of adversarial noise $\Delta$ and the smoothing parameter $h$ remain the same. The total number of oracle calls for $B > N$ changes to

\begin{align*}
 T &= \mathcal{O}\left(\max \left(\frac{LR^2}{\varepsilon}, \frac{\kappa\kappa'(p, d)d\sigma_*^2R^2}{\varepsilon^{2}}, ~\frac{\kappa_{\beta}^{2/(\beta - 1)}\kappa'(p, d)\kappa d^{2}\Delta^2L^{2/(\beta - 1)}R^{2 + 2/(\beta - 1)}}{\varepsilon^{2 + 2/(\beta - 1)}}\right)\right) \\ & = \mathcal{O}\left(\max \left(\frac{LR^2}{\varepsilon},~ \frac{\kappa\min \{q, \ln d\}d^{2 - 2/p}\sigma_*^2R^2}{\varepsilon^{2}},\frac{\min \{q, \ln d\}\kappa_{\beta}^{2/(\beta - 1)}\kappa d^{3 - 2/p}\Delta^2L^{2/(\beta - 1)}R^{2 + 2/(\beta - 1)}}{\varepsilon^{2 + 2/(\beta - 1)}}\right)\right).
\end{align*}
\end{remark}
 This remark  means that oracle complexity can be improved. For instance, optimizing in the $\ell_1$-norm can reduce the number of iterations by a factor of ${d}/{\ln{d}}$, which can be a significant speedup for large $d$.
The proof of this remark is given in Appendix \ref{app:rem_convergence_deterministic}.

\subsection{Stochastic adversarial noise}
% Write theorem with proof in the case of stochastic adversarial noise

Let us suppose that   a zero-order oracle is corrupted by a stochastic  noise $\delta$.
%\begin{equation}\label{eq:zero_order_stoch}
%f_\delta(x, \xi) = f(x, \xi) + \delta,
%\end{equation}
At each iteration we request two function evaluations on the inputs $x_t$ and $x_t^{'}$
\begin{equation}\label{eq:zero_order_stoch}
f_\delta(x_t, \xi) = f(x_t, \xi) + \delta_t \quad \mbox{and}\quad f_\delta(x_t^{'}, \xi) = f(x_t^{'}, \xi) + \delta_t^{'}.
\end{equation}
Then the gradient  can be approximated as follows
\begin{equation}\label{eq:grad_stoch_noise}
\mathbf{g}(x, \mathbf{e}, r, \xi) = \frac{d}{2h}(f(x + hr\mathbf{e}, \xi) - f(x - hr\mathbf{e}, \xi) + (\delta - \delta^{'}))K(r)\mathbf{e}.
\end{equation}
%We state the following assumption \cite{akhavan2022gradient} regarding stochastic adversarial noise.

\begin{assum}[Stochastic adversarial noise  \cite{akhavan2022gradient} ]
\label{ass:stochastic_noise}
For all $t = 1,2, ...$, it holds
\[\mathbb{E}[\delta_{t}^2] \le \Delta^2, \quad \mbox{and} \quad \mathbb{E}[{\delta}_{t}^{'2}] \le \Delta^2,\] ${\delta}_{t}$, ${\delta}_{t}^{'}$ are independent. 
\end{assum}
%of $\mathbf{e} \in S_2^d(1)$ and $r$ for kernel approximation

The next theorem states the convergence guarantee for AZO-SGD-HS with zero-order oracle \eqref{eq:zero_order_stoch}. We
 will use the shorthand notations $\kappa_\beta = \int |u|^\beta |K(u)|du$ and $\kappa = \int |K(u)^2|du$.

\begin{teo}
\label{teo:convergence_stochastic}
Let \(f\) satisfy Assumption \ref{ass:High-smoothness} with parameter $\beta$, the adversarial noise is stochastic and follows Assumption \ref{ass:stochastic_noise}. Let smoothing parameter be $h \le \left({\varepsilon}/{(\kappa_{\beta}LR})\right)^{{1}/{(\beta - 1)}}$ then    Assumption \ref{assum:overparam} is satisfied with gradient approximation  \eqref{eq_batched_grad}. Let the maximum admissible level of stochastic adversarial noise be 
\begin{itemize}
    \item
    $B \le N$
\begin{equation*}
\Delta \le \frac{\sigma_*\varepsilon^{{1}/{(\beta - 1)}}}{\kappa_{\beta}^{{1}/{(\beta - 1)}}\sqrt{d}R^{{1}/{(\beta - 1)}}L^{{1}/{(\beta - 1)}}}
\end{equation*}
    \item 
     $B > N$
\begin{equation*}
\Delta \le \frac{\varepsilon^{3/4 +{1}/{(\beta - 1)} }\sqrt{B}L^{1/4 - {1}/{(\beta - 1)}}}{\sqrt{\kappa\kappa'(p, d)}d\kappa_{\beta}^{{1}/{(\beta - 1)}}R^{1/2 + {1}/{(\beta - 1)}}}.
\end{equation*}
\end{itemize}
Let $x_N^{ag}$ be the output of the AZO-SGD-HS then 
\[\mathbb{E}[f(x_N^{ag}) - f^*] \le \varepsilon\] 
in at most
 the following number of iterations $N$ and  oracle calls $T$,
 where $R$ will denote a radius of a solution $x^*$: 
$\|x^* - x_0\| \le R$, and $x_0$ is a starting point:
 \begin{itemize}
     \item 

     $B = 1$
\begin{equation*}
N = \mathcal{O}\left(\max\left(\frac{LR^2}{\varepsilon}, \frac{d\kappa\sigma_*^2R^2}{\varepsilon^2}\right)\right), T = \mathcal{O}\left(\max\left(\frac{LR^2}{\varepsilon}, \frac{d\kappa\sigma_*^2R^2}{\varepsilon^2}\right)\right)
\end{equation*}

     \item 

$N > B > 1$

\begin{equation*}
N = \mathcal{O}\left(\max\left(\frac{LR^2}{B\varepsilon}, \sqrt{\frac{LR^2}{\varepsilon}}, \frac{d\kappa\sigma_*^2R^2}{B\varepsilon^2}\right)\right)
\end{equation*}
\begin{equation*}
T = \mathcal{O}\left(\max\left(\frac{LR^2}{\varepsilon}, B\sqrt{\frac{LR^2}{\varepsilon}}, \frac{d\kappa\sigma_*^2R^2}{\varepsilon^2}\right)\right) 
\end{equation*}

     \item 
$B = N$
\begin{equation*}
N = \mathcal{O} \left(\max \left(\sqrt{\frac{LR^2}{\varepsilon}}, \frac{\sqrt{d\kappa}\sigma_*R}{\varepsilon}\right)\right),T = \mathcal{O}\left(\max \left(\frac{LR^2}{\varepsilon}, \frac{d\kappa\sigma_*^2R^2}{\varepsilon^2}\right)\right) 
\end{equation*}

     \item 
     $B > N$

\begin{equation*}
N = \mathcal{O} \left(\sqrt{\frac{LR^2}{\varepsilon}}\right), T = \mathcal{O}\left(\max\left(\frac{LR^2}{\varepsilon}, \frac{d\kappa\sigma_*^2R^2}{\varepsilon^{2}}, \frac{\kappa_{\beta}^{2/(\beta - 1)}\kappa d^2\Delta^2L^{2/(\beta - 1)}R^{2 + 2/(\beta - 1)}}{\varepsilon^{2 + 2/(\beta - 1)}}\right)\right).
\end{equation*}
 \end{itemize}
\end{teo}

\textit{Sketch of the proof.}

\textbf{ Step 1} (bias of gradient approximation) 

In this step, we calculate the bias of the gradient approximation $\|\mathbb{E}[g(x_k, \xi, \mathbf{e})] - \nabla f(x_k)\|$. The main difference between them is 
\begin{align*}
b &= \|\mathbb{E}[\mathbf{g}(x_k, \mathbf{e}, r, \xi)] - \nabla f(x_k)\| = \left\|\mathbb{E}\left[d \frac{f_\delta(x_k + hr\mathbf{e}, \xi) - f_\delta(x_k - hr\mathbf{e}, \xi)}{2h} K(r) \mathbf{e}\right] - \nabla f(x_k)\right\| \\
& \stackrel{\text{\circlednum{1}}}{=}  \left\|\mathbb{E}\left[d \frac{f(x_k + hr\mathbf{e}, \xi)}{h}K(r) \mathbf{e}\right] + \mathbb{E}\left[\frac{\delta - \delta^{'}}{2h}K(r) \mathbf{e}\right] - \nabla f(x_k)\right\| \\
& \stackrel{\text{\circlednum{2}}}{=}\|\mathbb{E}\left[d \frac{f(x_k + hr\mathbf{e}, \xi)}{h}K(r) \mathbf{e}\right] - \nabla f(x_k)\|\stackrel{\text{\circlednum{3}}}{=}\|\mathbb{E}([\nabla f(x_k + h r \mathbf{u}, \xi)r K(r)] - \nabla f(x_k))\|\\
& \stackrel{\text{\circlednum{4}}}{\le}\sup_{z \in S_2^d(1)}\mathbb{E}[(\nabla_z f(x_k + h r \mathbf{u}, \xi) - \nabla_z f(x_k))r K(r)]  \\
& \stackrel{\text{\circlednum{5}}}{\le} \kappa_{\beta}h^{\beta - 1} \frac{L}{(l-1)!}\mathbb{E}[\|\mathbf{u}\|^{\beta - 1}] \le \kappa_{\beta}h^{\beta - 1} \frac{L}{(l-1)!}\frac{d}{d + \beta - 1}\lesssim \kappa_{\beta}L h^{\beta - 1},
\end{align*}
where $u \in B_2^d(1)$, \circlednum{1}: distribution of $\mathbf{e}$ is symmetric; \circlednum{2}: independence of the stochastic noise; \circlednum{3}: a version of the Stokes’
theorem \cite{zorich2016mathematical} (see Section 13.3.5, Exercise 14a); \circlednum{4}: norm of the gradient is the supremum of directional derivatives $\nabla_z f(x) = \lim_{\varepsilon \to 0}\frac{f(x + \varepsilon z) - f(x)}{\varepsilon}$;  \circlednum{5}: Taylor expansion.

\textbf{Step 2} (bounding the second moment of gradient approximation)

In this step, we calculate the second moment of the gradient approximation $\mathbb{E} [\|g(x^*, \xi, \mathbf{e})\|^2]$. The calculations are done similarly to the deterministic case.

\begin{align*}
\zeta^2 &= \mathbb{E} [\|\mathbf{g}(x^*, \mathbf{e}, r, \xi)\|^2] = \mathbb{E} \left [\left\|d \frac{f_\delta(x^* + hr\mathbf{e}, \xi) - f_\delta(x^* - hr\mathbf{e}, \xi)}{2h} K(r) \mathbf{e}\right\|^2\right]\\
&=\frac{d^2}{4h^2}\mathbb{E} [((f_\delta(x^* + hr\mathbf{e}, \xi) - f_\delta(x^* - hr\mathbf{e}, \xi)) K(r))^2] \\
& =\frac{d^2}{4h^2}\mathbb{E} [(f(x^* + hr\mathbf{e}, \xi)- f(x^* - hr\mathbf{e}, \xi) + (\delta - \delta^{'}))^2 (K(r))^2] \\
& \stackrel{\text{\circlednum{1}}}{\le} \frac{d^2 \kappa}{2h^2} (\mathbb{E} (f(x^* + hr\mathbf{e}, \xi)- f(x^* - hr\mathbf{e}, \xi))^2 + 2\Delta^2 )  \\
& \stackrel{\text{\circlednum{2}}}{\le} \frac{d^2 \kappa}{2h^2} (\frac{h^2}{d}\mathbb{E} (\|\nabla f(x^* + hr\mathbf{e}, \xi)+\nabla f(x^* - hr\mathbf{e}, \xi)\|^2) + 2\Delta^2) \\
& = \frac{d^2\kappa}{2h^2} (\frac{h^2}{d}\mathbb{E} (\|\nabla f(x^* + hr\mathbf{e}, \xi)+\nabla f(x^* - hr\mathbf{e}, \xi) \pm 2\nabla f(x^*, \xi)\|^2) + 2\Delta^2) \\
& \stackrel{\text{\circlednum{3}}}{\le} 4d\kappa\|\nabla f(x^*, \xi)\|^2 + 4d\kappa L^2h^2\mathbb{E}[\|\mathbf{e}\|^2] + \frac{d^2 \kappa\Delta^2}{h^2} \\
& \stackrel{\text{\circlednum{4}}}{\le} 4d\kappa\sigma_*^2 + 4d\kappa L^2h^2\mathbb{E}[\|\mathbf{e}\|^2] + \frac{d^2 \kappa\Delta^2}{h^2} = 4d\kappa\sigma_*^2 + 4d\kappa L^2h^2 + \frac{d^2\kappa\Delta^2}{h^2},
\end{align*}
\circlednum{1}: inequality of squared norm of the sum, inequality between positive random variables and independence of the noise; \circlednum{2}: Wirtinger-Poincare inequality; \circlednum{3}: L-smoothness function; \circlednum{4}: overparametrization assumption.

\qed

For the complete proof we refer to  Appendix \ref{app:teo_convergence_stochastic}.

Compared to Theorem \ref{teo:convergence_deterministic}, the Theorem \ref{teo:convergence_stochastic} shows similar results on iterative complexity $N$, total oracle calls $T$, and smoothing parameter $h$. The estimation on the maximum allowable level of adversarial noise in this case explicitly depends on batch size $B$, which demonstrates the effect of overbatching that already occurred for similar problems as in \cite{lobanov2024black}. The bound allows to significantly increase the quality of the solution with the increase in batch size.

\begin{remark}
\label{rem:convergence_stochastic}
Theorem \ref{teo:convergence_stochastic}  was proven for the Euclidean case but it cab generalized to the $\ell_p$-norm. In this case,  the  number of iterations $N$ and smoothing parameter $h$ remain the same. The total number of oracle calls changes to
\begin{align*}
 T &= \mathcal{O}\left(\max\left(\frac{LR^2}{\varepsilon},~ \frac{d\kappa\kappa'(p, d)\sigma_*^2R^2}{\varepsilon^{2}},~\frac{\kappa_{\beta}^{2/(\beta - 1)}\kappa\kappa'(p, d) d^2\Delta^2L^{2/(\beta - 1)}R^{2 + 2/(\beta - 1)}}{\varepsilon^{2 + 2/(\beta - 1)}}\right)\right) \\ & = \mathcal{O}\left(\max\left(\frac{LR^2}{\varepsilon}, ~\frac{\kappa\min \{q, \ln d\}d^{2 - 2/p}\sigma_*^2R^2}{\varepsilon^{2}},~\frac{\kappa_{\beta}^{2/(\beta - 1)}\kappa \min \{q, \ln d\}d^{3 - 2/p}\Delta^2L^{2/(\beta - 1)}R^{2 + 2/(\beta - 1)}}{\varepsilon^{2 + 2/(\beta - 1)}}\right)\right)
\end{align*}
and for $B > N$, the maximum admissible level of adversarial noise $\Delta$ changes to
\begin{equation*}
\Delta \le
\frac{\varepsilon^{3/4 +{1}/{(\beta - 1)} }\sqrt{B}L^{1/4 - {1}/{(\beta - 1)}}}{\sqrt{\kappa\min \{q, \ln d\}}d^{3/2 - 1/p}\kappa_{\beta}^{{1}/{(\beta - 1)}}R^{1/2 + {1}/{(\beta - 1)}}}.
\end{equation*}

\end{remark}
The proof of the this remark is given in Appendix \ref{app:rem_convergence_stochastic}. The remark also makes it possible  to reduce the total number of oracle calls in the $\ell_1$-norm. However, this reduces the maximum adversarial noise level that can potentially worsen the results.
In contrast, using the $\ell_\infty$ norm increases the total number of oracle calls and the maximum admissible level of adversarial noise, resulting in a trade-off between these parameters.

\section{Experiments}
\label{sec:experiments}

We evaluate the effectiveness of AZO-SGD-HS  by comparing it with with AZO-SGD   \cite{lobanov2023accelerated}, which does not utilize the property of high-order smoothness, on the logistic regression 

\begin{figure}[ht!]
    \centering
    \begin{subfigure}
        \centering
        \captionsetup{justification=centering}
        \caption*{Stochastic adversarial noise}
        \centering\includegraphics{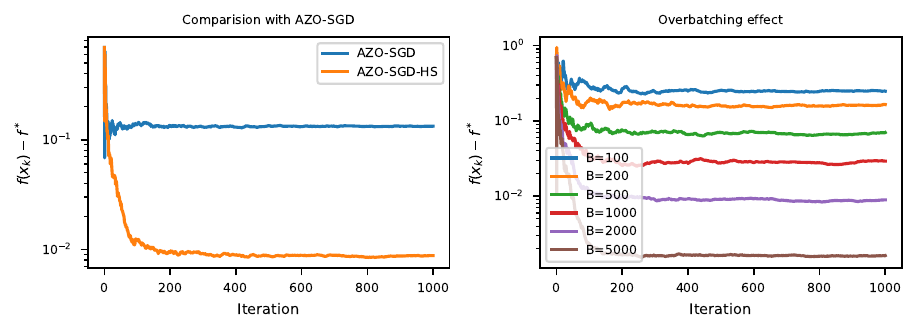} 
        % \captionsetup{justification=centering}
        % \caption*{Dual Extrapolation }
        \label{fig:fig2}
    \end{subfigure} 
        \begin{subfigure}
        \centering
          \captionsetup{justification=centering}
        \caption*{Deterministic adversarial noise}
        \centering\includegraphics{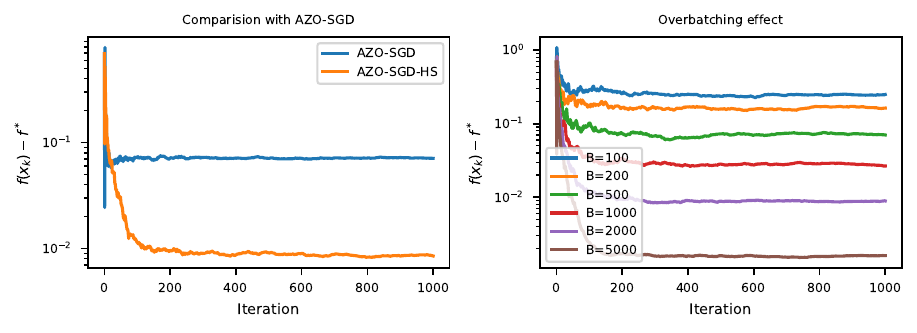}
        \label{fig:fig5}
    \end{subfigure}
    \caption{ Comparison of the proposed algorithm AZO-SGD-HS with AZO-SGD from \cite{lobanov2023accelerated} }\label{fig:deterministic_experiment}
\end{figure}

% Update experiments on Sirius and try ones with a more applied approach
%To support our theoretical results from Section \ref{sec:main} we evaluate the effectiveness of AZO-SGD-HS  on the logistic regression
%, a frequently used machine learning task at which our assumptions are fulfilled. The task can be formulated as the following optimization problem:
\begin{equation}\label{eq:logreg}
\min\limits_{w \in \R^{d+1}}  \sum_{i=1}^n -(y_i\log(p_i) + (1 - y_i)\log(1 - p_i)),
\end{equation}
where $p_i = 1/(1 + e^{\langle w, X_i\rangle})$, $X_i$ is the $i$-th example from the training dataset, and $y_i$ is the correct label of $X_i$, $n$  is the number of examples from training dataset, $w$ is the vector of parameters. The problem is convex w.r.t. $w$. Problem \eqref{eq:logreg} is smooth with parameter $L = {\lambda_{\max}(X^\top X)}/{(4n)}$, where $\lambda_{\max}(X^\top X)$ is the maximum eigenvalue of $X^\top X$. The high-order smoothness condition is also satisfied. The results are presented in  Figure \ref{fig:deterministic_experiment}. 

\paragraph{Synthetic data.} For binary classification, we
sampled $n=100$ examples with $d=1000$ features. The setup $d \gg n$ enjoys overparametrization. The dataset was sampled with the \textit{sklearn.datasets.make\_classification} procedure.

%We explored the dependence on the batch size $B$ value and compared our method for $\beta=3,4$ with AZO-SGD   \cite{bib:LobanovAccelerated}. 

\paragraph{Setup.} The parameters of the algorithms were chosen as follows: number of iterations $N = 1000$, adversarial noise level $\Delta = 0.0001$, $\beta=4$ for AZO-SGD-HS and theoretically calculated $L = 4.3$ for AZO-SGD \cite{lobanov2023accelerated}.  We  compare the methods  under batch size $B = 2000$. For demonstration of overbatching effect we considered $B = \{100, 200, 500, 1000, 2000, 5000\}$. To model deterministic adversarial noise we decreased the weights by a noisy normed gradient of the loss function w.r.t weights. 
%We explored the dependence on the batch size $B$ value and compared our method for  
To model stochastic adversarial noise, we used normal distribution  clipped to $[-\Delta, \Delta]$  with variance equal to $\Delta^2$. 
The kernel approximation was based on the Legendre polynomials \cite{bach2016highly}
\[K(r) = \sum_{m=0}^\beta p_m^{'}(0)p_m(r), \quad p_m(u) = \sqrt{2m+1}L_m(u),\] 
where $L_m$ is the $m$-th Legendre polynomial. 
For $\beta=4$
\[  K(r) =  \frac{15r}{4}(5-7r^2).\]
%For example, the following kernels for the first $\beta$ values are the following:
%\begin{equation}
%K(r) = \begin{cases}
%3r, \beta=1,2\\
%\frac{15r}{4}(5-7r^2), \beta=3,4\\
%\frac{105r}{64}(99r^4 - 126r^2+35), \beta=5,6
%\end{cases}
%\end{equation}

%We demonstrated an ``overbatching'' effect of the proposed method that occurs in better solution quality with the growth of batch size.
%The results of the experiments can be seen in Figure \ref{fig:stochastic_experiment}. They show similar effects, like in adversarial noise: an increase in the method's performance with the growth of the $B$ and superiority over AZO-SGD.

\section{Conclusion}
\label{sec:conclusion}

In this paper, we proposed a new gradient-free  algorithm (AZO-SGD-HS)  utilizing  a higher-order smoothness condition and designed for the overparametrized regime. The convergence of the proposed method was studied in the Euclidean case and then was generalized  to the non-Euclidean setup. The algorithm achieves optimal oracle call complexity, when the batch size $B$ is large enough. 
The convergence guarantee as well as the superiority of the proposed method over AZO-SGD \cite{lobanov2023accelerated} were confirmed by numerical experiments.

\bibliographystyle{unsrt}  
{\normalsize \bibliography{template}}

\begin{thebibliography}{10}

\bibitem{polyak1987introduction}
Boris Polyak.
\newblock {\em Introduction to Optimization}.
\newblock New York, Optimization Software, 1987.

\bibitem{granichin2003randomized}
ON~Granichin and BT~Polyak.
\newblock Randomized algorithms of an estimation and optimization under almost arbitrary noises.
\newblock {\em M.: Nauka}, 2003.

\bibitem{risteski2016algorithms}
Andrej Risteski and Yuanzhi Li.
\newblock Algorithms and matching lower bounds for approximately-convex optimization.
\newblock {\em Advances in Neural Information Processing Systems}, 29:4745--4753, 2016.

\bibitem{vasin2021stopping}
Artem Vasin, Alexander Gasnikov, and Vladimir Spokoiny.
\newblock Stopping rules for accelerated gradient methods with additive noise in gradient, 2021.

\bibitem{spall2003introduction}
James~C. Spall.
\newblock {\em Introduction to Stochastic Search and Optimization}.
\newblock John Wiley \& Sons, Inc., New York, NY, USA, 1 edition, 2003.

\bibitem{conn2009introduction}
A.~Conn, K.~Scheinberg, and L.~Vicente.
\newblock {\em Introduction to Derivative-Free Optimization}.
\newblock Society for Industrial and Applied Mathematics, 2009.

\bibitem{duchi2015optimal}
John~C. Duchi, Michael~I. Jordan, Martin~J. Wainwright, and Andre Wibisono.
\newblock Optimal rates for zero-order convex optimization: The power of two function evaluations.
\newblock {\em {IEEE} Trans. Information Theory}, 61(5):2788--2806, 2015.

\bibitem{shamir2017optimal}
Ohad Shamir.
\newblock An optimal algorithm for bandit and zero-order convex optimization with two-point feedback.
\newblock {\em Journal of Machine Learning Research}, 18:52:1--52:11, 2017.
\newblock First appeared in arXiv:1507.08752.

\bibitem{nesterov2017random}
Yurii Nesterov and Vladimir Spokoiny.
\newblock Random gradient-free minimization of convex functions.
\newblock {\em Found. Comput. Math.}, 17(2):527--566, April 2017.
\newblock First appeared in 2011 as CORE discussion paper 2011/16.

\bibitem{gasnikov2017stochastic}
A.~V. Gasnikov, E.~A. Krymova, A.~A. Lagunovskaya, I.~N. Usmanova, and F.~A. Fedorenko.
\newblock Stochastic online optimization. single-point and multi-point non-linear multi-armed bandits. convex and strongly-convex case.
\newblock {\em Automation and Remote Control}, 78(2):224--234, Feb 2017.
\newblock arXiv:1509.01679.

\bibitem{beznosikov2020gradient}
Aleksandr Beznosikov, Abdurakhmon Sadiev, and Alexander Gasnikov.
\newblock Gradient-free methods with inexact oracle for convex-concave stochastic saddle-point problem.
\newblock In {\em International Conference on Mathematical Optimization Theory and Operations Research}, pages 105--119. Springer, 2020.

\bibitem{gasnikov2022power}
Alexander Gasnikov, Anton Novitskii, Vasilii Novitskii, Farshed Abdukhakimov, Dmitry Kamzolov, Aleksandr Beznosikov, Martin Takac, Pavel Dvurechensky, and Bin Gu.
\newblock The power of first-order smooth optimization for black-box non-smooth problems.
\newblock {\em arXiv preprint arXiv:2201.12289}, 2022.

\bibitem{flaxman2004online}
Abraham~D Flaxman, Adam~Tauman Kalai, and H~Brendan McMahan.
\newblock Online convex optimization in the bandit setting: gradient descent without a gradient.
\newblock {\em arXiv preprint cs/0408007}, 2004.

\bibitem{bartlett2008high}
Peter Bartlett, Varsha Dani, Thomas Hayes, Sham Kakade, Alexander Rakhlin, and Ambuj Tewari.
\newblock High-probability regret bounds for bandit online linear optimization.
\newblock In {\em Proceedings of the 21st Annual Conference on Learning Theory-COLT 2008}, pages 335--342. Omnipress, 2008.

\bibitem{bubeck2012regret}
S\'ebastien Bubeck and Nicol\'o Cesa-Bianchi.
\newblock Regret analysis of stochastic and nonstochastic multi-armed bandit problems.
\newblock {\em Foundations and Trends® in Machine Learning}, 5(1):1--122, 2012.

\bibitem{ge2021anti}
Yunjie Ge, Qian Wang, Baolin Zheng, Xinlu Zhuang, Qi~Li, Chao Shen, and Cong Wang.
\newblock Anti-distillation backdoor attacks: Backdoors can really survive in knowledge distillation.
\newblock In {\em Proceedings of the 29th ACM International Conference on Multimedia}, pages 826--834, 2021.

\bibitem{dvinskikh2022gradient}
Darina Dvinskikh, Vladislav Tominin, Yaroslav Tominin, and Alexander Gasnikov.
\newblock Gradient-free optimization for non-smooth saddle point problems under adversarial noise.
\newblock {\em arXiv preprint arXiv:2202.06114}, 2022.

\bibitem{lobanov2023stochastic}
Aleksandr Lobanov.
\newblock Stochastic adversarial noise in the “black box” optimization problem.
\newblock In {\em International Conference on Optimization and Applications}, pages 60--71. Springer, 2023.

\bibitem{bayandina2018gradient}
Anastasia~Sergeevna Bayandina, Alexander~V Gasnikov, and Anastasia~A Lagunovskaya.
\newblock Gradient-free two-point methods for solving stochastic nonsmooth convex optimization problems with small non-random noises.
\newblock {\em Automation and Remote Control}, 79(8):1399--1408, 2018.

\bibitem{jacot2018neural}
Arthur Jacot, Franck Gabriel, and Cl{\'e}ment Hongler.
\newblock Neural tangent kernel: Convergence and generalization in neural networks.
\newblock {\em Advances in neural information processing systems}, 31, 2018.

\bibitem{allen2019learning}
Zeyuan Allen-Zhu, Yuanzhi Li, and Yingyu Liang.
\newblock Learning and generalization in overparameterized neural networks, going beyond two layers.
\newblock {\em Advances in neural information processing systems}, 32, 2019.

\bibitem{belkin2019reconciling}
Mikhail Belkin, Daniel Hsu, Siyuan Ma, and Soumik Mandal.
\newblock Reconciling modern machine-learning practice and the classical bias--variance trade-off.
\newblock {\em Proceedings of the National Academy of Sciences}, 116(32):15849--15854, 2019.

\bibitem{hoffmann2022training}
Jordan Hoffmann, Sebastian Borgeaud, Arthur Mensch, Elena Buchatskaya, Trevor Cai, Eliza Rutherford, Diego de~Las Casas, Lisa~Anne Hendricks, Johannes Welbl, Aidan Clark, et~al.
\newblock Training compute-optimal large language models.
\newblock {\em arXiv preprint arXiv:2203.15556}, 2022.

\bibitem{eldan2023tinystories}
Ronen Eldan and Yuanzhi Li.
\newblock Tinystories: How small can language models be and still speak coherent english?
\newblock {\em arXiv preprint arXiv:2305.07759}, 2023.

\bibitem{lobanov2023accelerated}
Aleksandr Lobanov and Alexander Gasnikov.
\newblock Accelerated zero-order sgd method for solving the black box optimization problem under “overparametrization” condition.
\newblock In {\em International Conference on Optimization and Applications}, pages 72--83. Springer, 2023.

\bibitem{polyak1990optimal}
Boris~Teodorovich Polyak and Aleksandr~Borisovich Tsybakov.
\newblock Optimal order of accuracy of search algorithms in stochastic optimization.
\newblock {\em Problemy Peredachi Informatsii}, 26(2):45--53, 1990.

\bibitem{bach2016highly}
Francis Bach and Vianney Perchet.
\newblock Highly-smooth zero-th order online optimization.
\newblock In {\em Conference on Learning Theory}, pages 257--283. PMLR, 2016.

\bibitem{akhavan2023gradient}
Arya Akhavan, Evgenii Chzhen, Massimiliano Pontil, and Alexandre~B Tsybakov.
\newblock Gradient-free optimization of highly smooth functions: improved analysis and a new algorithm.
\newblock {\em arXiv preprint arXiv:2306.02159}, 2023.

\bibitem{woodworth2021even}
Blake~E Woodworth and Nathan Srebro.
\newblock An even more optimal stochastic optimization algorithm: minibatching and interpolation learning.
\newblock {\em Advances in Neural Information Processing Systems}, 34:7333--7345, 2021.

\bibitem{zorich2016mathematical}
Vladimir~Antonovich Zorich and Octavio Paniagua.
\newblock {\em Mathematical analysis II}, volume 220.
\newblock Springer, 2016.

\bibitem{akhavan2022gradient}
Arya Akhavan, Evgenii Chzhen, Massimiliano Pontil, and Alexandre Tsybakov.
\newblock A gradient estimator via l1-randomization for online zero-order optimization with two point feedback.
\newblock {\em Advances in Neural Information Processing Systems}, 35:7685--7696, 2022.

\bibitem{lobanov2024black}
Aleksandr Lobanov, Nail Bashirov, and Alexander Gasnikov.
\newblock The “black-box” optimization problem: Zero-order accelerated stochastic method via kernel approximation.
\newblock {\em Journal of Optimization Theory and Applications}, pages 1--36, 2024.

\end{thebibliography}


\appendix
\section{Appendix}
\subsection{Proof of Theorem \ref{teo:convergence_deterministic}}
\label{app:teo_convergence_deterministic} 
We reuse the bias and the second moment of gradient approximation from the main part of the paper. They are:
\begin{align*}
& b \lesssim \kappa_{\beta}L h^{\beta - 1} + \frac{d\Delta}{h} \\
& \zeta^2 \le 4d\kappa\sigma_*^2 + 4d\kappa L^2h^2 + \frac{d^2\kappa\Delta^2}{h^2}.
\end{align*}

\textbf{Convergence}

We use the bound for biased oracle from \cite{lobanov2023accelerated} and substitute the bias and the second moment found earlier.

\begin{align*}
& \mathbb{E} [f(x_N^{ag}) - f^*] = c \left( \frac{LR^2}{N^2} + \frac{LR^2}{BN} + \frac{\zeta R}{\sqrt{BN}} + b R + \frac{b^2}{2L}N\right) \\
& \lesssim \frac{LR^2}{N^2} + \frac{LR^2}{BN} + \frac{\left(\sqrt{d\kappa\sigma_*^2} + \sqrt{d\kappa L^2h^2} + \sqrt{\frac{d^2\kappa\Delta^2}{h^2}}\right)R}{\sqrt{BN}} + \left(\kappa_{\beta}L h^{\beta - 1} + \frac{d\Delta}{h}\right) R \\
& +\frac{\left(\kappa_{\beta}L h^{\beta - 1} + \frac{d\Delta}{h}\right)^2}{L}N \le\frac{LR^2}{N^2} + \frac{LR^2}{BN} + \frac{\left(\sqrt{d\kappa}\sigma_* + \sqrt{d\kappa}Lh + \sqrt{\kappa}\frac{d\Delta}{h}\right)R}{\sqrt{BN}} \\
& +\left(\kappa_{\beta}L h^{\beta - 1} +
\frac{d\Delta}{h}\right) R + \frac{\kappa_{\beta}^2L^2 h^{2(\beta - 1)} + \frac{d^2\Delta^2}{h^2}}{L}N \\
& =\underbrace{\frac{LR^2}{N^2}}_{\text{\circlednum{1}}} + \underbrace{\frac{LR^2}{BN}}_{\text{\circlednum{2}}} + \underbrace{\frac{\sqrt{d\kappa}\sigma_*R}{\sqrt{BN}}}_{\text{\circlednum{3}}} + \underbrace{\frac{\sqrt{d\kappa}LhR}{\sqrt{BN}}}_{\text{\circlednum{4}}} + \underbrace{\frac{\sqrt{\kappa}d\Delta R}{h\sqrt{BN}}}_{\text{\circlednum{5}}} + \underbrace{\kappa_{\beta}L h^{\beta - 1}R}_{\text{\circlednum{6}}} + \underbrace{\frac{d\Delta}{h}R}_{\text{\circlednum{7}}} + \underbrace{\kappa_{\beta}^2L h^{2(\beta - 1)}}_{\text{\circlednum{8}}} + \underbrace{\frac{d^2\Delta^2N}{h^2L}}_{\text{\circlednum{9}}}.
\end{align*}

We consider several cases for the $B$ values and find all parameter values for all cases 

\framebox{Case 1: $B=1$} \par

\begin{align*}
& \mathbb{E} [f(x_N^{ag}) - f^*] = \frac{LR^2}{N^2} + \frac{LR^2}{N} + \frac{\sqrt{d\kappa}\sigma_*R}{\sqrt{N}} + \frac{\sqrt{d\kappa}LhR}{\sqrt{N}} + \frac{\sqrt{\kappa}d\Delta R}{h\sqrt{N}} + \kappa_{\beta}L h^{\beta - 1}R\\
& + \frac{d\Delta}{h}R + \kappa_{\beta}^2L h^{2(\beta - 1)} + \frac{d^2\Delta^2N}{h^2L}\le 9\varepsilon.
\end{align*}

We bound all terms by $\varepsilon$ and find bounds for needed parameters.

From \circlednum{1}, \circlednum{2} and \circlednum{3}, we find $N$

\begin{align*}
& \frac{LR^2}{N^2} \le \varepsilon \implies N \ge  \sqrt{\frac{LR^2}{\varepsilon}} = \mathcal{O}\left( \sqrt{\frac{LR^2}{\varepsilon}}\right) \\
& \frac{LR^2}{N} \le \varepsilon \implies N \ge \frac{LR^2}{\varepsilon} \\
& \frac{\sqrt{d\kappa}\sigma_*R}{\sqrt{N}}\le \varepsilon \implies N \ge \frac{d\kappa\sigma_*^2R^2}{\varepsilon^2}.
\end{align*}

Note that $\frac{LR^2}{\varepsilon} > \sqrt{\frac{LR^2}{\varepsilon}}$. That leads to the following expression for $N$
$$
N = \mathcal{O}\left(\max\left(\frac{LR^2}{\varepsilon}, \frac{d\kappa\sigma_*^2R^2}{\varepsilon^2}\right)\right).
$$

From \circlednum{4}, \circlednum{6} and \circlednum{8}, we find $h$

\begin{align*}
& \frac{\sqrt{d\kappa}LhR}{\sqrt{N}} \le \varepsilon \implies h \le \frac{\sqrt{N}\varepsilon}{\sqrt{d\kappa}LR} = \max \left(\frac{\sqrt{\frac{LR^2}{\varepsilon}}\varepsilon}{\sqrt{d\kappa}LR}, \frac{\sqrt{\frac{d\kappa\sigma_*^2R^2}{\varepsilon^2}}\varepsilon}{\sqrt{d\kappa}LR} \right) = \max \left(\frac{\varepsilon^{1/2}}{\sqrt{d\kappa L}}, \frac{\sigma_*}{L} \right) \\
& \kappa_{\beta}^2L h^{2(\beta - 1)} \le \varepsilon  \implies h \le \left(\frac{\varepsilon}{\kappa_{\beta}^2L}\right)^{\frac{1}{2(\beta - 1)}}\\
& \kappa_{\beta}L h^{\beta - 1}R  \le \varepsilon \implies h \le \left(\frac{\varepsilon}{\kappa_{\beta}LR}\right)^{\frac{1}{\beta - 1}}\\
& h = \min \left(\max \left(\frac{\varepsilon^{1/2}}{\sqrt{d\kappa L}}, \frac{\sigma_*}{L} \right), \left(\frac{\varepsilon}{\kappa_{\beta}^2L}\right)^{\frac{1}{2(\beta - 1)}}, \left(\frac{\varepsilon}{\kappa_{\beta}LR}\right)^{\frac{1}{\beta - 1}}\right) = \left(\frac{\varepsilon}{\kappa_{\beta}LR}\right)^{\frac{1}{\beta - 1}}.
\end{align*}

From \circlednum{5}, \circlednum{7} and \circlednum{9}, we find $\Delta$:

\begin{align*}
& \frac{\sqrt{\kappa}d\Delta R}{h\sqrt{N}} \le \varepsilon \implies \Delta \le \frac{\varepsilon h\sqrt{N}}{\sqrt{\kappa}dR} \le 
\max\left(\frac{\varepsilon \left(\frac{\varepsilon}{\kappa_{\beta}LR}\right)^{\frac{1}{\beta - 1}}\sqrt{\frac{LR^2}{\varepsilon}}}{\sqrt{\kappa}dR}, \frac{\varepsilon \left(\frac{\varepsilon}{\kappa_{\beta}LR}\right)^{\frac{1}{\beta - 1}}\sqrt{\frac{d\kappa\sigma_*^2R^2}{\varepsilon^2}}}{\sqrt{\kappa}dR}\right) \\
& =\max\left(\frac{\varepsilon^{1/2 + \frac{1}{\beta - 1}}L^{1/2 - \frac{1}{\beta - 1}}}{\sqrt{\kappa}d(\kappa_{\beta}R)^{\frac{1}{\beta - 1}}}, \frac{\varepsilon^{\frac{1}{\beta - 1}}\sigma_*}{\sqrt{d}\kappa_{\beta}^{\frac{1}{\beta - 1}}L^{\frac{1}{\beta - 1}}R^{\frac{1}{\beta - 1}}}\right) = \frac{\varepsilon^{\frac{1}{\beta - 1}}\sigma_*}{\sqrt{d}\kappa_{\beta}^{\frac{1}{\beta - 1}}L^{\frac{1}{\beta - 1}}R^{\frac{1}{\beta - 1}}}\\
& \frac{d\Delta}{h}R  \le \varepsilon \implies \Delta \le \frac{h\varepsilon}{dR} \le \frac{\left(\frac{\varepsilon}{\kappa_{\beta}LR}\right)^{\frac{1}{\beta - 1}}\varepsilon}{dR} = \frac{\varepsilon^{1 + \frac{1}{\beta - 1}}}{d\kappa_{\beta}^{\frac{1}{\beta - 1}}R^{1 + \frac{1}{\beta - 1}} L^{\frac{1}{\beta - 1}}} \\
& \frac{d^2\Delta^2N}{h^2L} \le \varepsilon \implies \Delta \le \frac{\sqrt{\varepsilon L}h}{\sqrt{N}d} \le \frac{\sqrt{\varepsilon L}\left(\frac{\varepsilon}{\kappa_{\beta}LR}\right)^{\frac{1}{\beta - 1}}}{\sqrt[4]{\frac{LR^2}{\varepsilon}}d} = \frac{\varepsilon^{1/2 + 1/4 + \frac{1}{\beta - 1}}L^{1/2 - \frac{1}{\beta - 1} - 1/4}}{dR^{1/2 + \frac{1}{\beta - 1}}\kappa_{\beta}^{\frac{1}{\beta - 1}}} \\
& =\frac{\varepsilon^{3/4 + \frac{1}{\beta - 1}}L^{1/4 - \frac{1}{\beta - 1}}}{dR^{1/2 + \frac{1}{\beta - 1}}\kappa_{\beta}^{\frac{1}{\beta - 1}}} \\
& \Delta = \min \left\{ 
\frac{\varepsilon^{\frac{1}{\beta - 1}}\sigma_*}{\sqrt{d}\kappa_{\beta}^{\frac{1}{\beta - 1}}L^{\frac{1}{\beta - 1}}R^{\frac{1}{\beta - 1}}}, \frac{\varepsilon^{1 + \frac{1}{\beta - 1}}}{d\kappa_{\beta}^{\frac{1}{\beta - 1}}R^{1 + \frac{1}{\beta - 1}} L^{\frac{1}{\beta - 1}}}, \frac{\varepsilon^{3/4 + \frac{1}{\beta - 1}}L^{1/4 - \frac{1}{\beta - 1}}}{dR^{1/2 + \frac{1}{\beta - 1}}\kappa_{\beta}^{\frac{1}{\beta - 1}}}\right\} \\
& =\frac{\varepsilon^{1 + \frac{1}{\beta - 1}}}{d\kappa_{\beta}^{\frac{1}{\beta - 1}}R^{1 + \frac{1}{\beta - 1}} L^{\frac{1}{\beta - 1}}}.
\end{align*}

Then, we can find total number of
oracle calls
$$
T = N\cdot B =\mathcal{O}\left(\max\left(\frac{LR^2}{\varepsilon}, \frac{d\kappa\sigma_*^2R^2}{\varepsilon^2}\right)\right).
$$

\framebox{Case 2: $N > B > 1$} \par

\begin{align*}
& \mathbb{E} [f(x_N^{ag}) - f^*] = \frac{LR^2}{N^2} + \frac{LR^2}{BN} + \frac{\sqrt{d\kappa}\sigma_*R}{\sqrt{BN}} + \frac{\sqrt{d\kappa}LhR}{\sqrt{BN}} + \frac{\sqrt{\kappa}d\Delta R}{h\sqrt{BN}} + \kappa_{\beta}L h^{\beta - 1}R \\
& +\frac{d\Delta}{h}R + \kappa_{\beta}^2L h^{2(\beta - 1)} + \frac{d^2\Delta^2N}{h^2L} \le 9\varepsilon.
\end{align*}

From \circlednum{1}, \circlednum{2}, \circlednum{3}, we find $N$

\begin{align*}
& \frac{LR^2}{N^2} \le \varepsilon \implies N \ge  \sqrt{\frac{LR^2}{\varepsilon}} \\
& \frac{LR^2}{BN}\le \varepsilon \implies N \ge \frac{LR^2}{B\varepsilon}\\
& \frac{\sqrt{d\kappa}\sigma_*R}{\sqrt{BN}}\le \varepsilon \implies N \ge \frac{d\kappa\sigma_*^2R^2}{B\varepsilon^2}\\
& N = \mathcal{O}\left(\max\left(\frac{LR^2}{B\varepsilon}, \sqrt{\frac{LR^2}{\varepsilon}}, \frac{d\kappa\sigma_*^2R^2}{B\varepsilon^2}\right)\right).
\end{align*}

From \circlednum{4}, \circlednum{6}, \circlednum{8}, we find $h$
\begin{align*}
& \frac{\sqrt{d\kappa}LhR}{\sqrt{BN}} \le \varepsilon \implies h \le \frac{\varepsilon\sqrt{BN}}{\sqrt{d\kappa}LR} \le \frac{\varepsilon\sqrt{BN}}{\sqrt{d\kappa}LR} \le \max(\frac{\varepsilon\sqrt{\frac{LR^2}{\varepsilon}}}{\sqrt{d\kappa}LR}, \frac{\varepsilon\sqrt{B\sqrt{\frac{LR^2}{\varepsilon}}}}{\sqrt{d\kappa}LR}, \frac{\varepsilon\sqrt{B\frac{d\kappa\sigma_*^2R^2}{B\varepsilon^2}}}{\sqrt{d\kappa}LR}) \\
& =\max\left(\sqrt{\frac{\varepsilon}{d\kappa L}}, \frac{\varepsilon^{3/4}\sqrt{B}}{\sqrt{d\kappa R}L^{3/4}}, \frac{\sigma_*}{L}\right)\\
& \kappa_{\beta}^2L h^{2(\beta - 1)} \le \varepsilon  \implies h \le \left(\frac{\varepsilon}{\kappa_{\beta}^2L}\right)^{\frac{1}{2(\beta - 1)}}\\
& \kappa_{\beta}L h^{\beta - 1}R  \le \varepsilon \implies h \le \left(\frac{\varepsilon}{\kappa_{\beta}LR}\right)^{\frac{1}{\beta - 1}}\\
& h = \min\left(\left(\frac{\varepsilon}{\kappa_{\beta}LR}\right)^{\frac{1}{\beta - 1}}, \left(\frac{\varepsilon}{\kappa_{\beta}^2L}\right)^{\frac{1}{2(\beta - 1)}}, \max\left(\sqrt{\frac{\varepsilon}{d\kappa L}}, \frac{\varepsilon^{3/4}\sqrt{B}}{\sqrt{d\kappa R}L^{3/4}}, \frac{\sigma_*}{L}\right)\right) = \left(\frac{\varepsilon}{\kappa_{\beta}LR}\right)^{\frac{1}{\beta - 1}}.
\end{align*}

From \circlednum{5}, \circlednum{7}, \circlednum{9}, we find $\Delta$

\begin{align*}
& \frac{\sqrt{\kappa}d\Delta R}{h\sqrt{BN}}\le \varepsilon \implies \Delta\le \frac{h\sqrt{BN}\varepsilon}{\sqrt{\kappa}d R} \le \max \left(\frac{h\left(\sqrt{\frac{LR^2}{\varepsilon}}\right)\varepsilon}{\sqrt{\kappa}d R}, \frac{h\left(\sqrt{B\sqrt{\frac{LR^2}{\varepsilon}}}\right)\varepsilon}{\sqrt{\kappa}d R}\right) \\
& =\max \left(\frac{h(\sqrt{L\varepsilon})}{\sqrt{\kappa}d}, \frac{h(\sqrt{B\sqrt{L}})\varepsilon^{3/4}}{\sqrt{\kappa}d R^{1/2}}\right) \\
& =\max \left(\frac{\left(\frac{\varepsilon}{\kappa_{\beta}LR}\right)^{\frac{1}{\beta - 1}}\left(\sqrt{L\varepsilon}\right)}{\sqrt{\kappa}d}, \frac{(\frac{\varepsilon}{\kappa_{\beta}LR})^{\frac{1}{\beta - 1}}(\sqrt{B\sqrt{L}})\varepsilon^{3/4}}{\sqrt{\kappa}d R^{1/2}}\right) \\
& =\max \left(\frac{\varepsilon^{1/2 + \frac{1}{\beta - 1}}L^{1/2 -\frac{1}{\beta - 1}}}{\sqrt{\kappa}d\kappa_{\beta}^{\frac{1}{\beta - 1}}R^{\frac{1}{\beta - 1}}}, \frac{L^{1/4 - \frac{1}{\beta - 1}}\sqrt{B}\varepsilon^{3/4 + \frac{1}{\beta - 1}}}{\sqrt{\kappa}d\kappa_{\beta}^{\frac{1}{\beta - 1}}R^{1/2 + \frac{1}{\beta - 1}}}\right) = \frac{\varepsilon^{1/2 + \frac{1}{\beta - 1}}L^{1/2 -\frac{1}{\beta - 1}}}{\sqrt{\kappa}d\kappa_{\beta}^{\frac{1}{\beta - 1}}R^{\frac{1}{\beta - 1}}}\\
& \frac{d\Delta}{h}R  \le \varepsilon \implies \Delta \le \frac{h\varepsilon}{dR} \le \frac{\left(\frac{\varepsilon}{\kappa_{\beta}LR}\right)^{\frac{1}{\beta - 1}}\varepsilon}{dR} = \frac{\varepsilon^{1 + \frac{1}{\beta - 1}}}{d\kappa_{\beta}^{\frac{1}{\beta - 1}}R^{1 + \frac{1}{\beta - 1}} L^{\frac{1}{\beta - 1}}}\\
& \frac{d^2\Delta^2N}{h^2L} \le \varepsilon \implies \Delta \le \frac{\sqrt{\varepsilon L}h}{\sqrt{N}d} \le \frac{\sqrt{\varepsilon L}\left(\frac{\varepsilon}{\kappa_{\beta}LR}\right)^{\frac{1}{\beta - 1}}}{\sqrt[4]{\frac{LR^2}{\varepsilon}}d} = \frac{\varepsilon^{1/2 + 1/4 + \frac{1}{\beta - 1}}L^{1/2 - \frac{1}{\beta - 1} - 1/4}}{dR^{1/2 + \frac{1}{\beta - 1}}\kappa_{\beta}^{\frac{1}{\beta - 1}}} \\
& =\frac{\varepsilon^{3/4 + \frac{1}{\beta - 1}}L^{1/4 - \frac{1}{\beta - 1}}}{dR^{1/2 + \frac{1}{\beta - 1}}\kappa_{\beta}^{\frac{1}{\beta - 1}}} \\
& \Delta = \min\left(\frac{\varepsilon^{1 + \frac{1}{\beta - 1}}}{d\kappa_{\beta}^{\frac{1}{\beta - 1}}R^{1 + \frac{1}{\beta - 1}} L^{\frac{1}{\beta - 1}}}, \frac{\varepsilon^{1/2 + \frac{1}{\beta - 1}}L^{1/2 -\frac{1}{\beta - 1}}}{\sqrt{\kappa}d\kappa_{\beta}^{\frac{1}{\beta - 1}}R^{\frac{1}{\beta - 1}}}, \frac{\varepsilon^{3/4 + \frac{1}{\beta - 1}}L^{1/4 - \frac{1}{\beta - 1}}}{dR^{1/2 + \frac{1}{\beta - 1}}\kappa_{\beta}^{\frac{1}{\beta - 1}}}\right) \\
& =\frac{\varepsilon^{1 + \frac{1}{\beta - 1}}}{d\kappa_{\beta}^{\frac{1}{\beta - 1}}R^{1 + \frac{1}{\beta - 1}} L^{\frac{1}{\beta - 1}}}\\
& T = N\cdot B =\mathcal{O}\left(\max \left(\frac{LR^2}{\varepsilon}, B\sqrt{\frac{LR^2}{\varepsilon}}, \frac{d\kappa\sigma_*^2R^2}{\varepsilon^2}\right)\right).
\end{align*}

\framebox{Case 3: $B = N$} \par

\begin{align*}
& \mathbb{E} [f(x_N^{ag}) - f^*] =\frac{LR^2}{N^2} + \frac{LR^2}{N^2} + \frac{\sqrt{d\kappa}\sigma_*R}{N} + \frac{\sqrt{d\kappa}LhR}{N} + \frac{\sqrt{\kappa}d\Delta R}{hN} + \kappa_{\beta}L h^{\beta - 1}R +\\
& +\frac{d\Delta}{h}R + \kappa_{\beta}^2L h^{2(\beta - 1)} + \frac{d^2\Delta^2N}{h^2L} \le 9\varepsilon.
\end{align*}

From \circlednum{1}, \circlednum{3}, we find $N$

\begin{align*}
& \frac{LR^2}{N^2} \le \varepsilon \implies N \ge  \sqrt{\frac{LR^2}{\varepsilon}} \\
& \frac{\sqrt{d\kappa}\sigma_*R}{N}\le \varepsilon \implies N \ge \frac{\sqrt{d\kappa}\sigma_*R}{\varepsilon} \\
& N = \max \left(\sqrt{\frac{LR^2}{\varepsilon}}, \frac{\sqrt{d\kappa}\sigma_*R}{\varepsilon}\right).
\end{align*}

From \circlednum{4}, \circlednum{6}, \circlednum{8}, we find $h$
\begin{align*}
& \frac{\sqrt{d\kappa}LhR}{N}\le \varepsilon \implies h \le \frac{\varepsilon N}{\sqrt{d\kappa}LR} = \max \left( \frac{\varepsilon}{\sqrt{d\kappa}LR}\sqrt{\frac{LR^2}{\varepsilon}}, \frac{\varepsilon}{\sqrt{d\kappa}LR}\frac{\sqrt{d\kappa}\sigma_*R}{\varepsilon}\right) \\
& =\max \left( \sqrt{\frac{\varepsilon}{d\kappa L}}, \frac{\sigma_*}{L}\right) \\
& \kappa_{\beta}^2L h^{2(\beta - 1)} \le \varepsilon  \implies h \le \left(\frac{\varepsilon}{\kappa_{\beta}^2L}\right)^{\frac{1}{2(\beta - 1)}} \\
& \kappa_{\beta}L h^{\beta - 1}R  \le \varepsilon \implies h \le \left(\frac{\varepsilon}{\kappa_{\beta}LR}\right)^{\frac{1}{\beta - 1}}\\
& h = \min \left(\max \left( \sqrt{\frac{\varepsilon}{d\kappa L}}, \frac{\sigma_*}{L}\right), \left(\frac{\varepsilon}{\kappa_{\beta}^2L}\right)^{\frac{1}{2(\beta - 1)}}, \left(\frac{\varepsilon}{\kappa_{\beta}LR}\right)^{\frac{1}{\beta - 1}}\right) = \left(\frac{\varepsilon}{\kappa_{\beta}LR}\right)^{\frac{1}{\beta - 1}}.
\end{align*}

From \circlednum{5}, \circlednum{7}, \circlednum{9}, we find $\Delta$

\begin{align*}
&\frac{\sqrt{\kappa}d\Delta R}{hN} \le \varepsilon \implies \Delta \le \frac{hN\varepsilon}{\sqrt{\kappa}d R} = \max \left(\frac{\left(\frac{\varepsilon}{\kappa_{\beta}LR}\right)^{\frac{1}{\beta - 1}}\varepsilon}{\sqrt{\kappa}d R}\sqrt{\frac{LR^2}{\varepsilon}}, \frac{\left(\frac{\varepsilon}{\kappa_{\beta}LR}\right)^{\frac{1}{\beta - 1}}\varepsilon}{\sqrt{\kappa}d R}\frac{\sqrt{d\kappa}\sigma_*R}{\varepsilon}\right) \\
& = \max \left(\frac{\varepsilon^{1/2 + \frac{1}{\beta - 1}}L^{1/2 - \frac{1}{\beta - 1}}}{\sqrt{\kappa}d\kappa_{\beta}^{\frac{1}{\beta - 1}}R^{\frac{1}{\beta - 1}}}, \frac{\sigma_*\varepsilon^{\frac{1}{\beta - 1}}}{\kappa_{\beta}^{\frac{1}{\beta - 1}}\sqrt{d}R^{\frac{1}{\beta - 1}}L^{\frac{1}{\beta - 1}}}\right) = \frac{\sigma_*\varepsilon^{\frac{1}{\beta - 1}}}{\kappa_{\beta}^{\frac{1}{\beta - 1}}\sqrt{d}R^{\frac{1}{\beta - 1}}L^{\frac{1}{\beta - 1}}}\\
& \frac{d\Delta}{h}R  \le \varepsilon \implies \Delta \le \frac{h\varepsilon}{dR} \le \frac{\left(\frac{\varepsilon}{\kappa_{\beta}LR}\right)^{\frac{1}{\beta - 1}}\varepsilon}{dR} = \frac{\varepsilon^{1 + \frac{1}{\beta - 1}}}{d\kappa_{\beta}^{\frac{1}{\beta - 1}}R^{1 + \frac{1}{\beta - 1}} L^{\frac{1}{\beta - 1}}}\\
& \frac{d^2\Delta^2N}{h^2L} \le \varepsilon \implies \Delta \le \frac{\sqrt{\varepsilon L}h}{\sqrt{N}d} \le \min \left(\frac{\sqrt{\varepsilon L}\left(\frac{\varepsilon}{\kappa_{\beta}LR}\right)^{\frac{1}{\beta - 1}}}{\sqrt[4]{\frac{LR^2}{\varepsilon}}d}, \frac{\sqrt{\varepsilon L}\left(\frac{\varepsilon}{\kappa_{\beta}LR}\right)^{\frac{1}{\beta - 1}}}{\sqrt{\frac{\sqrt{d\kappa}\sigma_*R}{\varepsilon}}d}
\right) \\
& =\min \left(\frac{\varepsilon^{1/2 + 1/4 + \frac{1}{\beta - 1}}L^{1/2 - \frac{1}{\beta - 1} - 1/4}}{dR^{1/2 + \frac{1}{\beta - 1}}\kappa_{\beta}^{\frac{1}{\beta - 1}}}, \frac{\varepsilon^{1 + \frac{1}{\beta - 1}}L^{1/2 - \frac{1}{\beta - 1}}}{d^{5/4}\kappa_{\beta}^{\frac{1}{\beta - 1}}\sqrt{\sqrt{\kappa}\sigma_*}R^{1/2 + \frac{1}{\beta - 1} }}\right) \\
& =\min \left(\frac{\varepsilon^{3/4 + \frac{1}{\beta - 1}}L^{1/4 - \frac{1}{\beta - 1}}}{dR^{1/2 + \frac{1}{\beta - 1}}\kappa_{\beta}^{\frac{1}{\beta - 1}}}, \frac{\varepsilon^{1 + \frac{1}{\beta - 1}}L^{1/2 - \frac{1}{\beta - 1}}}{d^{5/4}\kappa_{\beta}^{\frac{1}{\beta - 1}}\sqrt{\sqrt{\kappa}\sigma_*}R^{1/2 + \frac{1}{\beta - 1}}}\right) =  \frac{\varepsilon^{1 + \frac{1}{\beta - 1}}L^{1/2 - \frac{1}{\beta - 1}}}{d^{5/4}\kappa_{\beta}^{\frac{1}{\beta - 1}}\sqrt{\sqrt{\kappa}\sigma_*}R^{1/2 + \frac{1}{\beta - 1}}} \\
& \Delta = \min \left(\frac{\varepsilon^{1 + \frac{1}{\beta - 1}}L^{1/2 - \frac{1}{\beta - 1}}}{d^{5/4}\kappa_{\beta}^{\frac{1}{\beta - 1}}\sqrt{\sqrt{\kappa}\sigma_*}R^{1/2 + \frac{1}{\beta - 1}}}, \frac{\varepsilon^{1 + \frac{1}{\beta - 1}}}{d\kappa_{\beta}^{\frac{1}{\beta - 1}}R^{1 + \frac{1}{\beta - 1}} L^{\frac{1}{\beta - 1}}}, \frac{\sigma_*\varepsilon^{\frac{1}{\beta - 1}}}{\kappa_{\beta}^{\frac{1}{\beta - 1}}\sqrt{d}R^{\frac{1}{\beta - 1}}L^{\frac{1}{\beta - 1}}}\right) \\
& =\frac{\varepsilon^{1 + \frac{1}{\beta - 1}}}{d\kappa_{\beta}^{\frac{1}{\beta - 1}}R^{1 + \frac{1}{\beta - 1}} L^{\frac{1}{\beta - 1}}}\\
& T = N\cdot B =\mathcal{O}\left(\max \left(\frac{LR^2}{\varepsilon}, \frac{d\kappa\sigma_*^2R^2}{\varepsilon^2}\right)\right).
\end{align*}

\framebox{Case 4: $B > N$} \par

\begin{align*}
& \mathbb{E} [f(x_N^{ag}) - f^*] = \frac{LR^2}{N^2} + \frac{LR^2}{BN} + \frac{\sqrt{d\kappa}\sigma_*R}{\sqrt{BN}} + \frac{\sqrt{d\kappa}LhR}{\sqrt{BN}} + \frac{\sqrt{\kappa}d\Delta R}{h\sqrt{BN}} + \kappa_{\beta}L h^{\beta - 1}R \\
& +\frac{d\Delta}{h}R + \kappa_{\beta}^2L h^{2(\beta - 1)} + \frac{d^2\Delta^2N}{h^2L} \le 9\varepsilon.
\end{align*}

From \circlednum{1}, \circlednum{2}, \circlednum{3}, we find $N$ and $B$. Note that if $B > N$ then $\frac{LR^2}{N^2} > \frac{LR^2}{BN}$

\begin{align*}
& \frac{LR^2}{N^2} \le \varepsilon \implies B > N \ge  \sqrt{\frac{LR^2}{\varepsilon}}\\
& \frac{\sqrt{d\kappa}\sigma_*R}{\sqrt{BN}} \le \varepsilon \implies B \ge \frac{d\kappa\sigma_*^2R^2}{N\varepsilon^2} = \mathcal{O}\left(\frac{d\kappa\sigma_*^2R}{L^{1/2}\varepsilon^{3/2}}\right).
\end{align*}

From \circlednum{4}, \circlednum{6}, \circlednum{8}, we find $h$
\begin{align*}
& \kappa_{\beta}^2L h^{2(\beta - 1)} \le \varepsilon  \implies h \le \left(\frac{\varepsilon}{\kappa_{\beta}^2L}\right)^{\frac{1}{2(\beta - 1)}}\\
& \kappa_{\beta}L h^{\beta - 1}R  \le \varepsilon \implies h \le \left(\frac{\varepsilon}{\kappa_{\beta}LR}\right)^{\frac{1}{\beta - 1}}\\
& \frac{\sqrt{d\kappa}LhR}{\sqrt{BN}} \le \varepsilon \implies h \le \frac{\varepsilon\sqrt{BN}}{\sqrt{d\kappa}LR} \le \frac{\varepsilon\sqrt{BN}}{\sqrt{d\kappa}LR} \le \frac{\varepsilon\sqrt{B}\sqrt[4]{\frac{LR^2}{\varepsilon}}}{\sqrt{d\kappa}LR} = \frac{\varepsilon^{3/4}\sqrt{B}}{\sqrt{d\kappa} L^{3/4}R^{1/2}}\\
& h = \min\left(\left(\frac{\varepsilon}{\kappa_{\beta}^2L}\right)^{\frac{1}{2(\beta - 1)}}, \left(\frac{\varepsilon}{\kappa_{\beta}LR}\right)^{\frac{1}{\beta - 1}}, \frac{\varepsilon^{3/4}\sqrt{B}}{\sqrt{d\kappa} L^{3/4}R^{1/2}}\right) = \left(\frac{\varepsilon}{\kappa_{\beta}LR}\right)^{\frac{1}{\beta - 1}}.
\end{align*}

From \circlednum{5}, \circlednum{7}, \circlednum{9}, we find $\Delta$

\begin{align*}
& \frac{\sqrt{\kappa}d\Delta R}{h\sqrt{BN}}\le \varepsilon \implies \Delta\le \frac{h\sqrt{BN}\varepsilon}{\sqrt{\kappa}d R} \le \frac{\left(\frac{\varepsilon}{\kappa_{\beta}LR}\right)^{\frac{1}{\beta - 1}}\left(\sqrt{B}\sqrt[4]{\frac{LR^2}{\varepsilon}}\right)\varepsilon}{\sqrt{\kappa}d R} = \frac{\varepsilon^{3/4 + \frac{1}{\beta - 1}}\sqrt{B}L^{1/4 - \frac{1}{\beta - 1}}}{\sqrt{\kappa}d\kappa_{\beta}^{\frac{1}{\beta - 1}}R^{1/2 + \frac{1}{\beta - 1}}}\\
& \frac{d\Delta}{h}R  \le \varepsilon \implies \Delta \le \frac{h\varepsilon}{dR} \le \frac{\left(\frac{\varepsilon}{\kappa_{\beta}LR}\right)^{\frac{1}{\beta - 1}}\varepsilon}{dR} = \frac{\varepsilon^{1 + \frac{1}{\beta - 1}}}{d\kappa_{\beta}^{\frac{1}{\beta - 1}}R^{1 + \frac{1}{\beta - 1}} L^{\frac{1}{\beta - 1}}}\\
& \frac{d^2\Delta^2N}{h^2L} \le \varepsilon \implies \Delta \le \frac{\sqrt{\varepsilon L}h}{\sqrt{N}d} \le \frac{\sqrt{\varepsilon L}\left(\frac{\varepsilon}{\kappa_{\beta}LR}\right)^{\frac{1}{\beta - 1}}}{\sqrt[4]{\frac{LR^2}{\varepsilon}}d} = \frac{\varepsilon^{1/2 + 1/4 + \frac{1}{\beta - 1}}L^{1/2 - \frac{1}{\beta - 1} - 1/4}}{dR^{1/2 + \frac{1}{\beta - 1}}\kappa_{\beta}^{\frac{1}{\beta - 1}}} = \\
& =\frac{\varepsilon^{3/4 + \frac{1}{\beta - 1}}L^{1/4 - \frac{1}{\beta - 1}}}{dR^{1/2 + \frac{1}{\beta - 1}}\kappa_{\beta}^{\frac{1}{\beta - 1}}}\\
& \Delta = \min \left\{ 
\frac{\varepsilon^{3/4 + \frac{1}{\beta - 1}}\sqrt{B}L^{1/4 - \frac{1}{\beta - 1}}}{\sqrt{\kappa}d\kappa_{\beta}^{\frac{1}{\beta - 1}}R^{1/2 + \frac{1}{\beta - 1}}}, \frac{\varepsilon^{1 + \frac{1}{\beta - 1}}}{d\kappa_{\beta}^{\frac{1}{\beta - 1}}R^{1 + \frac{1}{\beta - 1}} L^{\frac{1}{\beta - 1}}}, \frac{\varepsilon^{3/4 + \frac{1}{\beta - 1}}L^{1/4 - \frac{1}{\beta - 1}}}{dR^{1/2 + \frac{1}{\beta - 1}}\kappa_{\beta}^{\frac{1}{\beta - 1}}}\right\} \\
& =\frac{\varepsilon^{1 + \frac{1}{\beta - 1}}}{d\kappa_{\beta}^{\frac{1}{\beta - 1}}R^{1 + \frac{1}{\beta - 1}} L^{\frac{1}{\beta - 1}}}.
\end{align*}

We can get the batch size $B$ via the maximum noise level $\Delta$ from \circlednum{5}:

$$
\frac{\sqrt{\kappa}d\Delta R}{h\sqrt{BN}}\le \varepsilon \implies B \ge \frac{\kappa d^2\Delta^2 R^2}{\varepsilon^2h^2N} \ge \frac{\kappa d^2\Delta^2 R^2}{\varepsilon^2\left(\frac{\varepsilon}{\kappa_{\beta}LR}\right)^{\frac{2}{\beta - 1}}\sqrt{\frac{LR^2}{\varepsilon}}} \ge \frac{\kappa_{\beta}^{\frac{2}{\beta - 1}}\kappa d^2\Delta^2R^{1 + \frac{2}{\beta - 1}}}{\varepsilon^{3/2 + \frac{2}{\beta - 1}}L^{1/2 - \frac{2}{\beta - 1}}}.
$$

That leads to the following $B$
\begin{align*}
& B = \mathcal{O}\left(\max\left(\sqrt{\frac{LR^2}{\varepsilon}}, \frac{d\kappa\sigma_*^2R}{L^{1/2}\varepsilon^{3/2}}, \frac{\kappa_{\beta}^{\frac{2}{\beta - 1}}\kappa d^2\Delta^2R^{1 + \frac{2}{\beta - 1}}}{\varepsilon^{3/2 + \frac{2}{\beta - 1}}L^{1/2 - \frac{2}{\beta - 1}}}\right)\right)\\
& T = B \cdot N = \mathcal{O}\left(\max \left(\frac{LR^2}{\varepsilon}, \frac{d\kappa\sigma_*^2R^2}{\varepsilon^{2}}, \frac{\kappa_{\beta}^{\frac{2}{\beta - 1}}\kappa d^2L^{\frac{2}{\beta - 1}}\Delta^2R^{2 + \frac{2}{\beta - 1}}}{\varepsilon^{2 + \frac{2}{\beta - 1}}}\right)\right).
\end{align*}

That ends the proof. \qed

\subsection{Proof of Remark \ref{rem:convergence_deterministic}}
\label{app:rem_convergence_deterministic}

\textbf{Bias of gradient approximation}

\begin{align*}
& b = \left\|\mathbb{E}\left[\mathbf{g}(x_k, \mathbf{e}, r, \xi)\right] - \nabla f(x_k)\right\|_p = \left\|\mathbb{E}\left[d \frac{f_\delta(x_k + hr\mathbf{e}, \xi) - f_\delta(x_k - hr\mathbf{e}, \xi)}{2h} K(r) \mathbf{e}\right] - \nabla f(x_k)\right\|_p \\
& \stackrel{\text{\circlednum{1}}}{=} \left\|\mathbb{E}\left[d \frac{f(x_k + hr\mathbf{e}, \xi)}{h}K(r) \mathbf{e} + \delta(x_k + h \mathbf{e})\right] - \nabla f(x_k)\right\|_p\\
& \stackrel{\text{\circlednum{2}}}{\le} \left\|\mathbb{E}\left[d \frac{f(x_k + hr\mathbf{e}, \xi)}{h}K(r) \mathbf{e}\right] - \nabla f(x_k)\right\|_p + \frac{d\Delta}{h}\\
& \stackrel{\text{\circlednum{3}}}{=}\|\mathbb{E}([\nabla f(x_k + h r \mathbf{u}, \xi)r K(r)] - \nabla f(x_k))\|_p + \frac{d\Delta}{h} \\
& \stackrel{\text{\circlednum{4}}}{=} \sup_{z \in S_p^d(1)}\mathbb{E}[(\nabla_z f(x_k + h r \mathbf{u}, \xi) - \nabla_z f(x_k))r K(r)] + \frac{d\Delta}{h} \\
& \stackrel{\text{\circlednum{5}}}{\le} \kappa_{\beta}h^{\beta - 1} \frac{L}{(l-1)!}\mathbb{E}[\|\mathbf{u}\|_p^{\beta - 1}] + \frac{d\Delta}{h} \lesssim \kappa_{\beta}L h^{\beta - 1} + \frac{d\Delta}{h},
\end{align*}
where $\mathbf{u} \in B_2^d(1)$, \circlednum{1}: distribution of $e$ is symmetric; \circlednum{2}: independence of the stochastic noise; \circlednum{3}: a version of the Stokes’
theorem \cite{zorich2016mathematical} (see Section 13.3.5, Exercise 14a); \circlednum{4}: norm of the gradient is the supremum of directional derivatives $\nabla_z f(x) = \lim_{\varepsilon \to 0}\frac{f(x + \varepsilon z) - f(x)}{\varepsilon}$;  \circlednum{5}: Taylor expansion. 

\textbf{Bounding second moment of gradient approximation}

Note that $\mathbb{E}[\|\mathbf{e}\|_p^2]\le \kappa'(p, d) = \min \{q, \ln d\}d^{2/q - 1} = \min \{q, \ln d\}d^{1 - 2/p}$ and 
\begin{align*}
& \zeta^2 = \mathbb{E}  [\|\mathbf{g}(x^*, \mathbf{e}, r, \xi)\|_p^2 ] = \mathbb{E} \left[\left\|d \frac{f_\delta(x^* + hr\mathbf{e}, \xi) - f_\delta(x^* - hr\mathbf{e}, \xi)}{2h} K(r) \mathbf{e}\right\|_p^2\right] \\
& =\frac{d^2}{4h^2}\mathbb{E}  [((f_\delta(x^* + hr\mathbf{e}, \xi) - f_\delta(x^* - hr\mathbf{e}, \xi)) K(r))^2\|\mathbf{e}\|_p^2] \\
& = \frac{d^2}{4h^2}\mathbb{E}  [(f(x^* + hr\mathbf{e}, \xi)- f(x^* - hr\mathbf{e}, \xi) + \delta(x^* + hr\mathbf{e}) - \delta(x^* - hr\mathbf{e}))^2 K^2(r)\|\mathbf{e}\|_p^2] \\
& \stackrel{\text{\circlednum{1}}}{\le}  \frac{d^2\kappa}{2h^2}\mathbb{E}  [(f(x^* + hr\mathbf{e}, \xi)- f(x^* - hr\mathbf{e}, \xi))^2\|\mathbf{e}\|_p^2] + \mathbb{E}  [(\delta(x^* + hr\mathbf{e}) - \delta(x^* - hr\mathbf{e}))^2\|\mathbf{e}\|_p^2]  \\
& \stackrel{\text{\circlednum{2}}}{\le} \frac{d^2\kappa}{2h^2}(\mathbb{E}  [(f(x^* + hr\mathbf{e}, \xi)- f(x^* - hr\mathbf{e}, \xi))^2\|\mathbf{e}\|_p^2] + 2\Delta^2 \mathbb{E}  [\|\mathbf{e}\|_p^2]) \\
& \le\frac{d^2\kappa\kappa'(p, d)}{2h^2}(\mathbb{E}  [(f(x^* + hr\mathbf{e}, \xi)- f(x^* - hr\mathbf{e}, \xi))^2] + 2\Delta^2 ) \\
& \stackrel{\text{\circlednum{3}}}{\le} \frac{d^2 \kappa\kappa'(p, d)}{2h^2} (\frac{h^2}{d}\mathbb{E} (\|\nabla f(x^* + hr\mathbf{e}, \xi)+\nabla f(x^* - hr\mathbf{e}, \xi)\|_2^2) + 2\Delta^2) \\
& = \frac{d^2\kappa\kappa'(p, d)}{2h^2} (\frac{h^2}{d}\mathbb{E} (\|\nabla f(x^* + hr\mathbf{e}, \xi)+\nabla f(x^* - hr\mathbf{e}, \xi) \pm 2\nabla f(x^*, \xi)\|_2^2) + 2\Delta^2)  \\
& \stackrel{\text{\circlednum{4}}}{\le} 4d\kappa\kappa'(p, d)\|\nabla f(x^*, \xi)\|_2^2 + 4d\kappa\kappa'(p, d) L^2h^2\mathbb{E}[\|\mathbf{e}\|_2^2] + \frac{d^2 \kappa\kappa'(p, d)\Delta^2}{h^2} \\
& \stackrel{\text{\circlednum{5}}}{\le} 4d\kappa\kappa'(p, d)\sigma_*^2 + 4d\kappa\kappa'(p, d) L^2h^2\mathbb{E}[\|\mathbf{e}\|_2^2] + \frac{d^2 \kappa\Delta^2}{h^2}\\
& =4d\kappa\kappa'(p, d)\sigma_*^2 + 4d\kappa\kappa'(p, d) L^2h^2 + \frac{d^2\kappa\kappa'(p, d)\Delta^2}{h^2},
\end{align*}
\circlednum{1}: inequality of squared norm of the sum and inequality between positive random variables; \circlednum{2}: bounded noise $|\delta| < \Delta$;  \circlednum{3}: Wirtinger-Poincare inequality; \circlednum{4}: L-smoothness function; \circlednum{5}: overparametrisation assumption. 

\textbf{Convergence}

\begin{align*}
& \mathbb{E} [f(x_N^{ag}) - f^*] = c \left( \frac{LR^2}{N^2} + \frac{LR^2}{BN} + \frac{\zeta R}{\sqrt{BN}} + b R + \frac{b^2}{2L}N\right) \\
& \lesssim \frac{LR^2}{N^2} + \frac{LR^2}{BN} + \frac{\left(\sqrt{d\kappa\kappa'(p, d)\sigma_*^2} + \sqrt{d\kappa \kappa'(p, d)L^2h^2} + \sqrt{\frac{d^2\kappa\kappa'(p, d)\Delta^2}{h^2}}\right)R}{\sqrt{BN}}\\
& +\left(\kappa_{\beta}L h^{\beta - 1} + \frac{d\Delta}{h}\right) R+\frac{(\kappa_{\beta}L h^{\beta - 1} + \frac{d\Delta}{h})^2}{L}N \le \frac{LR^2}{N^2} + \frac{LR^2}{BN}\\
& +\frac{\left(\sqrt{d\kappa\kappa'(p, d)}\sigma_* + \sqrt{d\kappa\kappa'(p, d)}Lh + \sqrt{\kappa\kappa'(p, d)}\frac{d\Delta}{h}\right)R}{\sqrt{BN}}
+\left(\kappa_{\beta}L h^{\beta - 1} + \frac{d\Delta}{h}\right) R \\
& +\frac{\kappa_{\beta}^2L^2 h^{2(\beta - 1)} + \frac{d^2\Delta^2}{h^2}}{L}N = \frac{LR^2}{N^2} + \frac{LR^2}{BN} + \frac{\sqrt{d\kappa\kappa'(p, d)}\sigma_*R}{\sqrt{BN}}\\
& +\frac{\sqrt{d\kappa\kappa'(p, d)}LhR}{\sqrt{BN}} +\frac{\sqrt{\kappa\kappa'(p, d)}d\Delta R}{h\sqrt{BN}} + \kappa_{\beta}L h^{\beta - 1}R + \frac{d\Delta}{h}R + \kappa_{\beta}^2L h^{2(\beta - 1)} + \frac{d^2\Delta^2N}{h^2L}.
\end{align*}

By performing the same transformations as in Theorem \ref{teo:convergence_deterministic}, we get the following parameters:

\framebox{Case 1: $B = 1$} \par

\begin{align*}
& N = \mathcal{O}\left(\max\left(\frac{LR^2}{\varepsilon}, \frac{d\kappa\kappa'(p, d)\sigma_*^2R^2}{\varepsilon^2}\right)\right)\\
& h = \left(\frac{\varepsilon}{\kappa_{\beta}LR}\right)^{\frac{1}{\beta - 1}}\\
& \Delta = \frac{\varepsilon^{1 + \frac{1}{\beta - 1}}}{d\kappa_{\beta}^{\frac{1}{\beta - 1}}R^{1 + \frac{1}{\beta - 1}} L^{\frac{1}{\beta - 1}}} \\
& T = \mathcal{O}\left(\max\left(\frac{LR^2}{\varepsilon}, \frac{d\kappa\kappa'(p, d)\sigma_*^2R^2}{\varepsilon^2}\right)\right).
\end{align*}

\framebox{Case 2: $N > B > 1$} \par
\begin{align*}
& N = \mathcal{O}\left(\max\left(\frac{LR^2}{B\varepsilon}, \sqrt{\frac{LR^2}{\varepsilon}}, \frac{d\kappa\kappa'(p, d)\sigma_*^2R^2}{B\varepsilon^2}\right)\right)\\
& h = \left(\frac{\varepsilon}{\kappa_{\beta}LR}\right)^{\frac{1}{\beta - 1}}\\
& \Delta = \frac{\varepsilon^{1 + \frac{1}{\beta - 1}}}{d\kappa_{\beta}^{\frac{1}{\beta - 1}}R^{1 + \frac{1}{\beta - 1}} L^{\frac{1}{\beta - 1}}}\\
& T = \mathcal{O}\left(\max\left(\frac{LR^2}{\varepsilon}, B\sqrt{\frac{LR^2}{\varepsilon}}, \frac{d\kappa\kappa'(p, d)\sigma_*^2R^2}{\varepsilon^2}\right)\right).
\end{align*}

\framebox{Case 3: $B = N$} \par

\begin{align*}
& N = \mathcal{O}\left(\max \left(\sqrt{\frac{LR^2}{\varepsilon}}, \frac{\sqrt{d\kappa\kappa'(p, d)}\sigma_*R}{\varepsilon}\right)\right)\\
& h = \left(\frac{\varepsilon}{\kappa_{\beta}LR}\right)^{\frac{1}{\beta - 1}}\\
& \Delta = \frac{\varepsilon^{1 + \frac{1}{\beta - 1}}}{d\kappa_{\beta}^{\frac{1}{\beta - 1}}R^{1 + \frac{1}{\beta - 1}} L^{\frac{1}{\beta - 1}}}\\
& T = \mathcal{O}\left(\max \left(\frac{LR^2}{\varepsilon}, \frac{d\kappa\kappa'(p, d)\sigma_*^2R^2}{\varepsilon^2}\right)\right).
\end{align*}

\framebox{Case 4: $B > N$} \par
\begin{align*}
& N = \mathcal{O}\left(\sqrt{\frac{LR^2}{\varepsilon}}\right)\\
& B = \mathcal{O}\left(\max\left(\sqrt{\frac{LR^2}{\varepsilon}}, \frac{d\kappa\kappa'(p, d)\sigma_*^2R}{\sqrt{L}\varepsilon^{3/2}}, \frac{\kappa_{\beta}^{\frac{2}{\beta - 1}}\kappa \kappa'(p, d)d^2\Delta^2R^{1 + \frac{2}{\beta - 1}}}{\varepsilon^{3/2 + \frac{2}{\beta - 1}}L^{1/2 - \frac{2}{\beta - 1}}}\right)\right)\\
& h =  \left(\frac{\varepsilon}{\kappa_{\beta}LR}\right)^{\frac{1}{\beta - 1}}\\
& \Delta = \frac{\varepsilon^{1 + \frac{1}{\beta - 1}}}{d\kappa_{\beta}^{\frac{1}{\beta - 1}}R^{1 + \frac{1}{\beta - 1}} L^{\frac{1}{\beta - 1}}}\\
& T = \mathcal{O}\left(\max \left(\frac{LR^2}{\varepsilon}, \frac{\kappa\kappa'(p, d)d\sigma_*^2R^2}{\varepsilon^{2}}, \frac{\kappa_{\beta}^{\frac{2}{\beta - 1}}\kappa'(p, d)\kappa d^{2}\Delta^2R^{1 + \frac{2}{\beta - 1}}}{\varepsilon^{3/2 + \frac{2}{\beta - 1}}L^{1/2 - \frac{2}{\beta - 1}}}\right)\right)=\\
& =\mathcal{O}\left(\max \left(\frac{LR^2}{\varepsilon}, \frac{\kappa\min \{q, \ln d\}d^{2 - 2/p}\sigma_*^2R^2}{\varepsilon^{2}}, \frac{\min \{q, \ln d\}\kappa_{\beta}^{\frac{2}{\beta - 1}}\kappa d^{3 - 2/p}\Delta^2R^{1 + \frac{2}{\beta - 1}}}{\varepsilon^{3/2 + \frac{2}{\beta - 1}}L^{1/2 - \frac{2}{\beta - 1}}}\right)\right).
\end{align*}

By substituting $\kappa'(p, d)$ into formulas, we prove the remark. \qed

\subsection{Proof of Theorem \ref{teo:convergence_stochastic}}
\label{app:teo_convergence_stochastic}

We reuse the bias and the second moment of gradient approximation from the main part of the paper. They are:
\begin{align*}
& b \lesssim \kappa_{\beta}L h^{\beta - 1}\\
& \zeta^2 \le 4d\kappa\sigma_*^2 + 4d\kappa L^2h^2 + \frac{d^2\kappa\Delta^2}{h^2}.
\end{align*}

\textbf{Convergence}

We reemploy the scheme for bounding $\mathbb{E} [f(x_N^{ag}) - f^*]$ from Theorem \ref{teo:convergence_deterministic}.

\begin{align*}
& \mathbb{E} [f(x_N^{ag}) - f^*] = c \left( \frac{LR^2}{N^2} + \frac{LR^2}{BN} + \frac{\zeta R}{\sqrt{BN}} + b R + \frac{b^2}{2L}N\right)  \\
& \lesssim \frac{LR^2}{N^2} + \frac{LR^2}{BN} + \frac{\left(\sqrt{d\kappa\sigma_*^2} + \sqrt{d\kappa L^2h^2} + \sqrt{\frac{d^2\kappa\Delta^2}{h^2}}\right)R}{\sqrt{BN}} + (\kappa_{\beta}L h^{\beta - 1}) R + \frac{(\kappa_{\beta}L h^{\beta - 1})^2}{L}N \\
& \le\frac{LR^2}{N^2} + \frac{LR^2}{BN} + \frac{\left(\sqrt{d\kappa}\sigma_* + \sqrt{d\kappa}Lh + \sqrt{\kappa}\frac{d\Delta}{h}\right)R}{\sqrt{BN}} + (\kappa_{\beta}L h^{\beta - 1}) R + \frac{\kappa_{\beta}^2L^2 h^{2(\beta - 1)}}{L}N \\
& =\underbrace{\frac{LR^2}{N^2}}_{\text{\circlednum{1}}} + \underbrace{\frac{LR^2}{BN}}_{\text{\circlednum{2}}} + \underbrace{\frac{\sqrt{d\kappa}\sigma_*R}{\sqrt{BN}}}_{\text{\circlednum{3}}} + \underbrace{\frac{\sqrt{d\kappa}LhR}{\sqrt{BN}}}_{\text{\circlednum{4}}} + \underbrace{\frac{\sqrt{\kappa}d\Delta R}{h\sqrt{BN}}}_{\text{\circlednum{5}}} + \underbrace{\kappa_{\beta}L h^{\beta - 1}R}_{\text{\circlednum{6}}} + \underbrace{\kappa_{\beta}^2L h^{2(\beta - 1)}}_{\text{\circlednum{7}}}.
\end{align*}

\framebox{Case 1: $B = 1$} \par

$$
\mathbb{E} [f(x_N^{ag}) - f^*]=\frac{LR^2}{N^2} + \frac{LR^2}{N} + \frac{\sqrt{d\kappa}\sigma_*R}{\sqrt{N}} + \frac{\sqrt{d\kappa}LhR}{\sqrt{N}} + \frac{\sqrt{\kappa}d\Delta R}{h\sqrt{N}} + \kappa_{\beta}L h^{\beta - 1}R + \kappa_{\beta}^2L h^{2(\beta - 1)} \le 7\varepsilon.
$$
From \circlednum{1}, \circlednum{2}, \circlednum{3}, we find $N$
\begin{align*}
& \frac{LR^2}{N^2} \le \varepsilon \implies N \ge  \sqrt{\frac{LR^2}{\varepsilon}} \\
& \frac{LR^2}{N} \le \varepsilon \implies N \ge \frac{LR^2}{\varepsilon}\\
& \frac{\sqrt{d\kappa}\sigma_*R}{\sqrt{N}}\le \varepsilon \implies N \ge \frac{d\kappa\sigma_*^2R^2}{\varepsilon^2}\\
& N = \mathcal{O}\left(\max\left(\frac{LR^2}{\varepsilon}, \frac{d\kappa\sigma_*^2R^2}{\varepsilon^2}\right)\right).
\end{align*}

From \circlednum{4}, \circlednum{6}, \circlednum{7}, we find $h$

\begin{align*}
& \frac{\sqrt{d\kappa}LhR}{\sqrt{N}} \le \varepsilon \implies h \le \frac{\sqrt{N}\varepsilon}{\sqrt{d\kappa}LR} \le \max \left(\frac{\sqrt{\frac{LR^2}{\varepsilon}}\varepsilon}{\sqrt{d\kappa}LR}, \frac{\sqrt{\frac{d\kappa\sigma_*^2R^2}{\varepsilon^2}}\varepsilon}{\sqrt{d\kappa}LR} \right) = \max \left(\frac{\varepsilon^{1/2}}{\sqrt{d\kappa L}}, \frac{\sigma_*}{L} \right) \\
& \kappa_{\beta}^2L h^{2(\beta - 1)} \le \varepsilon  \implies h \le \left(\frac{\varepsilon}{\kappa_{\beta}^2L}\right)^{\frac{1}{2(\beta - 1)}}\\
& \kappa_{\beta}L h^{\beta - 1}R  \le \varepsilon \implies h \le \left(\frac{\varepsilon}{\kappa_{\beta}LR}\right)^{\frac{1}{\beta - 1}}\\
& h = \min \left(\max \left(\frac{\varepsilon^{1/2}}{\sqrt{d\kappa L}}, \frac{\sigma_*}{L} \right),\left(\frac{\varepsilon}{\kappa_{\beta}^2L}\right)^{\frac{1}{2(\beta - 1)}},\left(\frac{\varepsilon}{\kappa_{\beta}LR}\right)^{\frac{1}{\beta - 1}}\right) = \left(\frac{\varepsilon}{\kappa_{\beta}LR}\right)^{\frac{1}{\beta - 1}}.
\end{align*}

From \circlednum{5}, we find $\Delta$

\begin{align*}
& \frac{\sqrt{\kappa}d\Delta R}{h\sqrt{N}} \le \varepsilon \implies \Delta \le \frac{\varepsilon h\sqrt{N}}{\sqrt{\kappa}dR} \le 
\max\left(\frac{\varepsilon \left(\frac{\varepsilon}{\kappa_{\beta}LR}\right)^{\frac{1}{\beta - 1}}\sqrt{\frac{LR^2}{\varepsilon}}}{\sqrt{\kappa}dR}, \frac{\varepsilon \left(\frac{\varepsilon}{\kappa_{\beta}LR}\right)^{\frac{1}{\beta - 1}}\sqrt{\frac{d\kappa\sigma_*^2R^2}{\varepsilon^2}}}{\sqrt{\kappa}dR}\right) \\
& =\max\left(\frac{\varepsilon^{1/2 + \frac{1}{\beta - 1}}L^{1/2 - \frac{1}{\beta - 1}}}{\sqrt{\kappa}d(\kappa_{\beta}R)^{\frac{1}{\beta - 1}}}, \frac{\varepsilon^{\frac{1}{\beta - 1}}\sigma_*}{\sqrt{d}\kappa_{\beta}^{\frac{1}{\beta - 1}}L^{\frac{1}{\beta - 1}}R^{\frac{1}{\beta - 1}}}\right) = \frac{\varepsilon^{\frac{1}{\beta - 1}}\sigma_*}{\sqrt{d}\kappa_{\beta}^{\frac{1}{\beta - 1}}L^{\frac{1}{\beta - 1}}R^{\frac{1}{\beta - 1}}}\\
& T = B \cdot N = \mathcal{O}\left(\max\left(\frac{LR^2}{\varepsilon}, \frac{d\kappa\sigma_*^2R^2}{\varepsilon^2}\right)\right).
\end{align*}

\framebox{Case 2: $N > B > 1$} \par

$$
\mathbb{E} [f(x_N^{ag}) - f^*]=\frac{LR^2}{N^2} + \frac{LR^2}{BN} + \frac{\sqrt{d\kappa}\sigma_*R}{\sqrt{BN}} + \frac{\sqrt{d\kappa}LhR}{\sqrt{BN}} + \frac{\sqrt{\kappa}d\Delta R}{h\sqrt{BN}} + \kappa_{\beta}L h^{\beta - 1}R + \kappa_{\beta}^2L h^{2(\beta - 1)} \le 7\varepsilon.
$$
From \circlednum{1}, \circlednum{2}, \circlednum{3}, we find $N$
\begin{align*}
& \frac{LR^2}{N^2} \le \varepsilon \implies N \ge  \sqrt{\frac{LR^2}{\varepsilon}}\\
& \frac{LR^2}{BN}\le \varepsilon \implies N \ge \frac{LR^2}{B\varepsilon}\\
& \frac{\sqrt{d\kappa}\sigma_*R}{\sqrt{BN}}\le \varepsilon \implies N \ge \frac{d\kappa\sigma_*^2R^2}{B\varepsilon^2}\\
& N = \mathcal{O}\left(\max\left(\frac{LR^2}{B\varepsilon}, \sqrt{\frac{LR^2}{\varepsilon}}, \frac{d\kappa\sigma_*^2R^2}{B\varepsilon^2}\right)\right).
\end{align*}

From \circlednum{4}, \circlednum{6}, \circlednum{7}, we find $h$
\begin{align*}
& \frac{\sqrt{d\kappa}LhR}{\sqrt{BN}}\le \varepsilon \implies h \le \frac{\sqrt{BN} \varepsilon}{\sqrt{d\kappa}LR} = \max\left(\frac{\sqrt{B\frac{LR^2}{B\varepsilon}} \varepsilon}{\sqrt{d\kappa}LR}, \frac{\sqrt{B\sqrt{\frac{LR^2}{\varepsilon}}} \varepsilon}{\sqrt{d\kappa}LR}, \frac{\sqrt{B\frac{d\kappa\sigma_*^2R^2}{B\varepsilon^2}} \varepsilon}{\sqrt{d\kappa}LR}\right)\\
& =\max\left(\sqrt{\frac{\varepsilon}{d\kappa L}}, \frac{\varepsilon^{3/4}B^{1/2}}{\sqrt{d\kappa}L^{3/4}R^{1/2}}, \frac{\sigma_*}{L}\right)\\
& \kappa_{\beta}^2L h^{2(\beta - 1)} \le \varepsilon  \implies h \le \left(\frac{\varepsilon}{\kappa_{\beta}^2L}\right)^{\frac{1}{2(\beta - 1)}}\\
& \kappa_{\beta}L h^{\beta - 1}R  \le \varepsilon \implies h \le \left(\frac{\varepsilon}{\kappa_{\beta}LR}\right)^{\frac{1}{\beta - 1}}\\
& h = \min\left(\max\left(\sqrt{\frac{\varepsilon}{d\kappa L}}, \frac{\varepsilon^{3/4}B^{1/2}}{\sqrt{d\kappa}L^{3/4}R^{1/2}}, \frac{\sigma_*}{L}\right), \left(\frac{\varepsilon}{\kappa_{\beta}^2L}\right)^{\frac{1}{2(\beta - 1)}}, \left(\frac{\varepsilon}{\kappa_{\beta}LR}\right)^{\frac{1}{\beta - 1}}\right) = \left(\frac{\varepsilon}{\kappa_{\beta}LR}\right)^{\frac{1}{\beta - 1}}.
\end{align*}

From \circlednum{5}, we find $\Delta$
\begin{align*}
& \frac{\sqrt{\kappa}d\Delta R}{h\sqrt{BN}}\le \varepsilon \implies \Delta\le \frac{h\sqrt{BN}\varepsilon}{\sqrt{\kappa}d R} \\
& \le\max \left(\frac{h\left(\sqrt{\frac{LR^2}{\varepsilon}}\right)\varepsilon}{\sqrt{\kappa}d R}, \frac{h\left(\sqrt{B\sqrt{\frac{LR^2}{\varepsilon}}}\right)\varepsilon}{\sqrt{\kappa}d R}, \frac{h\left(\sqrt{B\frac{d\kappa\sigma_*^2R^2}{B\varepsilon^2}}\right)\varepsilon}{\sqrt{\kappa}d R}\right) \\
& =\max \left(\frac{h\left(\sqrt{L\varepsilon}\right)}{\sqrt{\kappa}d}, \frac{h\left(\sqrt{B\sqrt{L}}\right)\varepsilon^{3/4}}{\sqrt{\kappa}d R^{1/2}}, \frac{\sigma_*h}{\sqrt{d}}\right) \\
& =\max \left(\frac{\left(\frac{\varepsilon}{\kappa_{\beta}LR}\right)^{\frac{1}{\beta - 1}}\left(\sqrt{L\varepsilon}\right)}{\sqrt{\kappa}d}, \frac{\left(\frac{\varepsilon}{\kappa_{\beta}LR}\right)^{\frac{1}{\beta - 1}}\left(\sqrt{B\sqrt{L}}\right)\varepsilon^{3/4}}{\sqrt{\kappa}d R^{1/2}}, \frac{\sigma_*}{\sqrt{d}}\left(\frac{\varepsilon}{\kappa_{\beta}LR}\right)^{\frac{1}{\beta - 1}}\right)\\
& = \max \left(\frac{\varepsilon^{1/2 + \frac{1}{\beta - 1}}L^{1/2 -\frac{1}{\beta - 1}}}{\sqrt{\kappa}d\kappa_{\beta}^{\frac{1}{\beta - 1}}R^{\frac{1}{\beta - 1}}}, \frac{L^{1/4 - \frac{1}{\beta - 1}}\sqrt{B}\varepsilon^{3/4 + \frac{1}{\beta - 1}}}{\sqrt{\kappa}d\kappa_{\beta}^{\frac{1}{\beta - 1}}R^{1/2 + \frac{1}{\beta - 1}}}, \frac{\sigma_*}{\sqrt{d}}\left(\frac{\varepsilon}{\kappa_{\beta}LR}\right)^{\frac{1}{\beta - 1}}\right) = \frac{\sigma_*}{\sqrt{d}}\left(\frac{\varepsilon}{\kappa_{\beta}LR}\right)^{\frac{1}{\beta - 1}}\\
& T = N\cdot B =\mathcal{O}\left(\max\left(\frac{LR^2}{\varepsilon}, B\sqrt{\frac{LR^2}{\varepsilon}}, \frac{d\kappa\sigma_*^2R^2}{\varepsilon^2}\right)\right).
\end{align*}

\framebox{Case 3: $B = N$} \par

$$
\mathbb{E} [f(x_N^{ag}) - f^*]=\frac{LR^2}{N^2} + \frac{LR^2}{N^2} + \frac{\sqrt{d\kappa}\sigma_*R}{N} + \frac{\sqrt{d\kappa}LhR}{N} + \frac{\sqrt{\kappa}d\Delta R}{hN} + \kappa_{\beta}L h^{\beta - 1}R + \kappa_{\beta}^2L h^{2(\beta - 1)} \le 7\varepsilon
$$
From \circlednum{1}, \circlednum{2}, \circlednum{3}, we find $N$
\begin{align*}
& \frac{LR^2}{N^2}\le \varepsilon \implies N \ge \sqrt{\frac{LR^2}{\varepsilon}}\\
& \frac{\sqrt{d\kappa}\sigma_*R}{N} \le \varepsilon \implies N \ge \frac{\sqrt{d\kappa}\sigma_*R}{\varepsilon}\\
& N = \mathcal{O}\left(\max \left(\sqrt{\frac{LR^2}{\varepsilon}}, \frac{\sqrt{d\kappa}\sigma_*R}{\varepsilon}\right)\right).
\end{align*}

From \circlednum{4}, \circlednum{6}, \circlednum{7}, we find $h$
\begin{align*}
& \frac{\sqrt{d\kappa}LhR}{N}\le \varepsilon \implies h \le \frac{\varepsilon N}{\sqrt{d\kappa}LR} = \max \left( \frac{\varepsilon}{\sqrt{d\kappa}LR}\sqrt{\frac{LR^2}{\varepsilon}}, \frac{\varepsilon}{\sqrt{d\kappa}LR}\frac{\sqrt{d\kappa}\sigma_*R}{\varepsilon}\right) \\
& =\max \left( \sqrt{\frac{\varepsilon}{d\kappa L}}, \frac{\sigma_*}{L}\right)\\
& \kappa_{\beta}^2L h^{2(\beta - 1)} \le \varepsilon  \implies h \le \left(\frac{\varepsilon}{\kappa_{\beta}^2L}\right)^{\frac{1}{2(\beta - 1)}}\\
& \kappa_{\beta}L h^{\beta - 1}R  \le \varepsilon \implies h \le \left(\frac{\varepsilon}{\kappa_{\beta}LR}\right)^{\frac{1}{\beta - 1}}\\
& h = \min \left(\max \left( \sqrt{\frac{\varepsilon}{d\kappa L}}, \frac{\sigma_*}{L}\right), \left(\frac{\varepsilon}{\kappa_{\beta}^2L}\right)^{\frac{1}{2(\beta - 1)}}, \left(\frac{\varepsilon}{\kappa_{\beta}LR}\right)^{\frac{1}{\beta - 1}}\right) = \left(\frac{\varepsilon}{\kappa_{\beta}LR}\right)^{\frac{1}{\beta - 1}}.
\end{align*}

From \circlednum{5}, we find $\Delta$

\begin{align*}
& \frac{\sqrt{\kappa}d\Delta R}{hN} \le \varepsilon \implies \Delta \le \frac{hN\varepsilon}{\sqrt{\kappa}d R} = \max \left(\frac{\left(\frac{\varepsilon}{\kappa_{\beta}LR}\right)^{\frac{1}{\beta - 1}}\varepsilon}{\sqrt{\kappa}d R}\sqrt{\frac{LR^2}{\varepsilon}}, \frac{\left(\frac{\varepsilon}{\kappa_{\beta}LR}\right)^{\frac{1}{\beta - 1}}\varepsilon}{\sqrt{\kappa}d R}\frac{\sqrt{d\kappa}\sigma_*R}{\varepsilon}\right) \\
& = \max \left(\frac{\varepsilon^{1/2 + \frac{1}{\beta - 1}}L^{1/2 - \frac{1}{\beta - 1}}}{\sqrt{\kappa}d\kappa_{\beta}^{\frac{1}{\beta - 1}}R^{\frac{1}{\beta - 1}}}, \frac{\sigma_*\varepsilon^{\frac{1}{\beta - 1}}}{\kappa_{\beta}^{\frac{1}{\beta - 1}}\sqrt{d}R^{\frac{1}{\beta - 1}}L^{\frac{1}{\beta - 1}}}\right) = \frac{\sigma_*\varepsilon^{\frac{1}{\beta - 1}}}{\kappa_{\beta}^{\frac{1}{\beta - 1}}\sqrt{d}R^{\frac{1}{\beta - 1}}L^{\frac{1}{\beta - 1}}}\\
& T = N\cdot B =\mathcal{O}\left(\max \left(\frac{LR^2}{\varepsilon}, \frac{d\kappa\sigma_*^2R^2}{\varepsilon^2}\right)\right).
\end{align*}

\framebox{Case 4: $B > N$} \par

$$
\mathbb{E} [f(x_N^{ag}) - f^*]=\frac{LR^2}{N^2} + \frac{LR^2}{BN} + \frac{\sqrt{d\kappa}\sigma_*R}{\sqrt{BN}} + \frac{\sqrt{d\kappa}LhR}{\sqrt{BN}} + \frac{\sqrt{\kappa}d\Delta R}{h\sqrt{BN}} + \kappa_{\beta}L h^{\beta - 1}R + \kappa_{\beta}^2L h^{2(\beta - 1)} \le 7\varepsilon.
$$

From \circlednum{1}, \circlednum{2}, \circlednum{3}, we find $N$ and $B$. Note that if $B > N$ then $\frac{LR^2}{N^2} > \frac{LR^2}{BN}$
\begin{align*}
& \frac{LR^2}{N^2}\le \varepsilon \implies B > N \ge \sqrt{\frac{LR^2}{\varepsilon}}\\
& \frac{\sqrt{d\kappa}\sigma_*R}{\sqrt{BN}} \le \varepsilon \implies B \ge \frac{d\kappa\sigma_*^2R^2}{N\varepsilon^2} = \mathcal{O}\left(\frac{d\kappa\sigma_*^2R}{L^{1/2}\varepsilon^{3/2}}\right).
\end{align*}

From \circlednum{4}, \circlednum{6}, \circlednum{7}, we find $h$

\begin{align*}
& \frac{\sqrt{d\kappa}LhR}{\sqrt{BN}} \le \varepsilon \implies h \le \frac{\varepsilon\sqrt{BN}}{\sqrt{d\kappa}LR} \le \frac{\varepsilon\sqrt{B}\sqrt[4]{\frac{LR^2}{\varepsilon}}}{\sqrt{d\kappa}LR} = \frac{\varepsilon^{3/4}\sqrt{B}}{\sqrt{d\kappa}L^{3/4}R^{1/2}}\\
& \kappa_{\beta}^2L h^{2(\beta - 1)} \le \varepsilon  \implies h \le \left(\frac{\varepsilon}{\kappa_{\beta}^2L}\right)^{\frac{1}{2(\beta - 1)}}\\
& \kappa_{\beta}L h^{\beta - 1}R  \le \varepsilon \implies h \le \left(\frac{\varepsilon}{\kappa_{\beta}LR}\right)^{\frac{1}{\beta - 1}}\\
& h = \min\left(\left(\frac{\varepsilon}{\kappa_{\beta}^2L}\right)^{\frac{1}{2(\beta - 1)}}, \left(\frac{\varepsilon}{\kappa_{\beta}LR}\right)^{\frac{1}{\beta - 1}}, \frac{\varepsilon^{3/4}\sqrt{B}}{\sqrt{d\kappa}L^{3/4}R^{1/2}}\right) = \left(\frac{\varepsilon}{\kappa_{\beta}LR}\right)^{\frac{1}{\beta - 1}}.
\end{align*}

From \circlednum{5}, we find $\Delta$
\begin{align*}
& \frac{\sqrt{\kappa}d\Delta R}{h\sqrt{BN}}\le \varepsilon \implies \Delta\le \frac{h\sqrt{BN}\varepsilon}{\sqrt{\kappa}d R} \le \frac{\left(\frac{\varepsilon}{\kappa_{\beta}LR}\right)^{\frac{1}{\beta - 1}}\left(\sqrt{B}\sqrt[4]{\frac{LR^2}{\varepsilon}}\right)\varepsilon}{\sqrt{\kappa}d R} \\
& =\frac{\left(\frac{\varepsilon}{\kappa_{\beta}LR}\right)^{\frac{1}{\beta - 1}}\left(\sqrt{B}\sqrt[4]{\frac{LR^2}{\varepsilon}}\right)\varepsilon}{\sqrt{\kappa}d R}=\frac{\varepsilon^{3/4 +\frac{1}{\beta - 1} }\sqrt{B}L^{1/4 - \frac{1}{\beta - 1}}}{\sqrt{\kappa}d\kappa_{\beta}^{\frac{1}{\beta - 1}}R^{1/2 + \frac{1}{\beta - 1}}}
\end{align*}

or we can get the batch size $B$ via the maximum noise level $\Delta$:
\begin{align*}
& \frac{\sqrt{\kappa}d\Delta R}{h\sqrt{BN}}\le \varepsilon \implies B \ge \frac{\kappa d^2\Delta^2 R^2}{\varepsilon^2h^2N} \ge \frac{\kappa d^2\Delta^2 R^2}{\varepsilon^2\left(\frac{\varepsilon}{\kappa_{\beta}LR}\right)^{\frac{2}{\beta - 1}}\sqrt{\frac{LR^2}{\varepsilon}}} \\
& =\frac{\kappa_{\beta}^{\frac{2}{\beta - 1}}\kappa d^2\Delta^2R^{1 + \frac{2}{\beta - 1}}}{\varepsilon^{3/2 + \frac{2}{\beta - 1}}L^{1/2 - \frac{2}{\beta - 1}}} \\
& B = \mathcal{O}\left(\max\left(\sqrt{\frac{LR^2}{\varepsilon}}, \frac{d\kappa\sigma_*^2R}{L^{1/2}\varepsilon^{3/2}}, \frac{\kappa_{\beta}^{\frac{2}{\beta - 1}}\kappa d^2\Delta^2R^{1 + \frac{2}{\beta - 1}}}{\varepsilon^{3/2 + \frac{2}{\beta - 1}}L^{1/2 - \frac{2}{\beta - 1}}}\right)\right)\\
& T = N \cdot B = \mathcal{O}\left(\max\left(\frac{LR^2}{\varepsilon}, \frac{d\kappa\sigma_*^2R^2}{\varepsilon^{2}}, \frac{\kappa_{\beta}^{\frac{2}{\beta - 1}}\kappa d^2\Delta^2L^{\frac{2}{\beta - 1}}R^{2 + \frac{2}{\beta - 1}}}{\varepsilon^{2 + \frac{2}{\beta - 1}}}\right)\right).
\end{align*}

That ends the proof. \qed

\subsection{Proof of Remark \ref{rem:convergence_stochastic}}
\label{app:rem_convergence_stochastic}

\textbf{Bias of gradient approximation}
\begin{align*}
& b = \|\mathbb{E}[\mathbf{g}(x_k, \mathbf{e}, r, \xi)] - \nabla f(x_k)\|_p = \left\|\mathbb{E}\left[d \frac{f_\delta(x_k + hr\mathbf{e}, \xi) - f_\delta(x_k - hr\mathbf{e}, \xi)}{2h} K(r) \mathbf{e}\right] - \nabla f(x_k)\right\|_p \\
& \stackrel{\text{\circlednum{1}}}{=} \left\|\mathbb{E}\left[d \frac{f(x_k + hr\mathbf{e}, \xi)}{h}K(r) \mathbf{e} + \delta(x_k + h \mathbf{e})\right] - \nabla f(x_k)\right\|_p\\
& \stackrel{\text{\circlednum{2}}}{\le} \left\|\mathbb{E}\left[d \frac{f(x_k + hr\mathbf{e}, \xi)}{h}K(r) \mathbf{e}\right] - \nabla f(x_k)\right\|_p\\
& \stackrel{\text{\circlednum{3}}}{=}\|\mathbb{E}([\nabla f(x_k + h r \mathbf{u}, \xi)r K(r)] - \nabla f(x_k))\|_p\\
& \stackrel{\text{\circlednum{4}}}{=} \sup_{z \in S_p^d(1)}\mathbb{E}[(\nabla_z f(x_k + h r \mathbf{u}, \xi) - \nabla_z f(x_k))r K(r)] \\
& \stackrel{\text{\circlednum{5}}}{\le} \kappa_{\beta}h^{\beta - 1} \frac{L}{(l-1)!}\mathbb{E}[\|\mathbf{u}\|_p^{\beta - 1}] \lesssim \kappa_{\beta}L h^{\beta - 1},
\end{align*}
where $\mathbf{u} \in B_2^d(1)$, \circlednum{1}: distribution of $e$ is symmetric; \circlednum{2}: independence of the stochastic noise; \circlednum{3}: a version of the Stokes’
theorem \cite{zorich2016mathematical} (see Section 13.3.5, Exercise 14a); \circlednum{4}: norm of the gradient is the supremum of directional derivatives $\nabla_z f(x) = \lim_{\varepsilon \to 0}\frac{f(x + \varepsilon z) - f(x)}{\varepsilon}$;  \circlednum{5}: Taylor expansion. 

\textbf{Bounding second moment of gradient approximation}

The calculations can be taken from Remark \ref{rem:convergence_deterministic} without changes. 

\textbf{Convergence}

\begin{align*}
& \mathbb{E} [f(x_N^{ag}) - f^*] = c \left( \frac{LR^2}{N^2} + \frac{LR^2}{BN} + \frac{\zeta R}{\sqrt{BN}} + b R + \frac{b^2}{2L}N\right) \\
& \lesssim \frac{LR^2}{N^2} + \frac{LR^2}{BN} + \frac{\left(\sqrt{d\kappa\kappa'(p, d)\sigma_*^2} + \sqrt{d\kappa \kappa'(p, d)L^2h^2} + \sqrt{\frac{d^2\kappa\kappa'(p, d)\Delta^2}{h^2}}\right)R}{\sqrt{BN}} +(\kappa_{\beta}L h^{\beta - 1}) R\\
& + \frac{(\kappa_{\beta}L h^{\beta - 1})^2}{L}N \le \frac{LR^2}{N^2} + \frac{LR^2}{BN} + \frac{\left(\sqrt{d\kappa\kappa'(p, d)}\sigma_* + \sqrt{d\kappa\kappa'(p, d)}Lh + \sqrt{\kappa\kappa'(p, d)}\frac{d\Delta}{h}\right)R}{\sqrt{BN}} \\
& +(\kappa_{\beta}L h^{\beta - 1}) R + \frac{\kappa_{\beta}^2L^2 h^{2(\beta - 1)}}{L}N = \frac{LR^2}{N^2} + \frac{LR^2}{BN} + \frac{\sqrt{d\kappa\kappa'(p, d)}\sigma_*R}{\sqrt{BN}} \\
& +\frac{\sqrt{d\kappa\kappa'(p, d)}LhR}{\sqrt{BN}} +\frac{\sqrt{\kappa\kappa'(p, d)}d\Delta R}{h\sqrt{BN}} + \kappa_{\beta}L h^{\beta - 1}R + \kappa_{\beta}^2L h^{2(\beta - 1)} .
\end{align*}

By performing the same transformations as in Theorem \ref{teo:convergence_stochastic}, we get the following parameters:

\framebox{Case 1: $B = 1$} \par

\begin{align*}
& N = \mathcal{O}\left(\max\left(\frac{LR^2}{\varepsilon}, \frac{d\kappa\kappa'(p, d)\sigma_*^2R^2}{\varepsilon^2}\right)\right)\\
& h \le \left(\frac{\varepsilon}{\kappa_{\beta}LR}\right)^{\frac{1}{\beta - 1}}\\
& \Delta \le \frac{\varepsilon^{\frac{1}{\beta - 1}}\sigma_*}{\sqrt{d}\kappa_{\beta}^{\frac{1}{\beta - 1}}L^{\frac{1}{\beta - 1}}R^{\frac{1}{\beta - 1}}}\\
& T = \mathcal{O}\left(\max\left(\frac{LR^2}{\varepsilon}, \frac{d\kappa\kappa'(p, d)\sigma_*^2R^2}{\varepsilon^2}\right)\right).
\end{align*}

\framebox{Case 2: $N > B > 1$} \par

\begin{align*}
& N = \mathcal{O}\left(\max\left(\frac{LR^2}{B\varepsilon}, \sqrt{\frac{LR^2}{\varepsilon}}, \frac{d\kappa\kappa'(p, d)\sigma_*^2R^2}{B\varepsilon^2}\right)\right)\\
& h \le \left(\frac{\varepsilon}{\kappa_{\beta}LR}\right)^{\frac{1}{\beta - 1}}\\
& \Delta \le \frac{\varepsilon^{\frac{1}{\beta - 1}}\sigma_*}{\sqrt{d}\kappa_{\beta}^{\frac{1}{\beta - 1}}L^{\frac{1}{\beta - 1}}R^{\frac{1}{\beta - 1}}}\\
& T = \mathcal{O}\left(\max\left(\frac{LR^2}{\varepsilon}, B\sqrt{\frac{LR^2}{\varepsilon}}, \frac{d\kappa\kappa'(p, d)\sigma_*^2R^2}{\varepsilon^2}\right)\right).
\end{align*}

\framebox{Case 3: $B = N$} \par

\begin{align*}
& N = \mathcal{O}\left(\max \left(\sqrt{\frac{LR^2}{\varepsilon}}, \frac{\sqrt{d\kappa\kappa'(p, d)}\sigma_*R}{\varepsilon}\right)\right)\\
& h \le \left(\frac{\varepsilon}{\kappa_{\beta}LR}\right)^{\frac{1}{\beta - 1}}\\
& \Delta \le \frac{\sigma_*\varepsilon^{\frac{1}{\beta - 1}}}{\kappa_{\beta}^{\frac{1}{\beta - 1}}\sqrt{d}R^{\frac{1}{\beta - 1}}L^{\frac{1}{\beta - 1}}}\\
& T = \mathcal{O}\left(\max \left(\frac{LR^2}{\varepsilon}, \frac{d\kappa\kappa'(p, d)\sigma_*^2R^2}{\varepsilon^2}\right)\right).
\end{align*}

\framebox{Case 4: $B > N$} \par

\begin{align*}
& N = \mathcal{O} \left(\sqrt{\frac{LR^2}{\varepsilon}}\right)\\
& B = \mathcal{O}\left(\max\left(\sqrt{\frac{LR^2}{\varepsilon}}, \frac{d\kappa\kappa'(p, d)\sigma_*^2R}{L^{1/2}\varepsilon^{3/2}}, \frac{\kappa_{\beta}^{\frac{2}{\beta - 1}}\kappa \kappa'(p, d)d^2\Delta^2R^{1 + \frac{2}{\beta - 1}}}{\varepsilon^{3/2 + \frac{2}{\beta - 1}}L^{1/2 - \frac{2}{\beta - 1}}}\right)\right)\\
& h \le \left(\frac{\varepsilon}{\kappa_{\beta}LR}\right)^{\frac{1}{\beta - 1}}\\
& \Delta \le \frac{\varepsilon^{3/4 +\frac{1}{\beta - 1} }\sqrt{B}L^{1/4 - \frac{1}{\beta - 1}}}{\sqrt{\kappa\kappa'(p, d)}d\kappa_{\beta}^{\frac{1}{\beta - 1}}R^{1/2 + \frac{1}{\beta - 1}}}\\
& T = \mathcal{O}\left(\max\left(\frac{LR^2}{\varepsilon}, \frac{d\kappa\kappa'(p, d)\sigma_*^2R^2}{\varepsilon^{2}}, \frac{\kappa_{\beta}^{\frac{2}{\beta - 1}}\kappa\kappa'(p, d) d^2\Delta^2L^{\frac{2}{\beta - 1}}R^{2 + \frac{2}{\beta - 1}}}{\varepsilon^{2 + \frac{2}{\beta - 1}}}\right)\right).
\end{align*}

By substituting $\kappa'(p, d)$ into formulas, we prove the remark. \qed
\end{document}